\documentclass[11pt,oneside]{article}
\usepackage[T1]{fontenc} 
\usepackage[nottoc]{tocbibind}
\usepackage[bottom=4cm]{geometry} % see geometry.pdf on how to lay out the page. There's lots.
\usepackage{float} %für grafik plazierung nötig

\usepackage{amsfonts}
\usepackage{epsfig}
\usepackage{url}
\usepackage{amssymb,latexsym,amsmath,amsthm,verbatim}
\usepackage{graphicx,epsfig,epstopdf,amssymb,color}
\usepackage{bm}

\newcommand{\veps}{\varepsilon}
\newcommand{\R}{\mathbb{R}} 
\newcommand{\C}{\mathbb{C}}
\newcommand{\K}{\mathbb{K}}

\newcommand{\ba}{\begin{array}}
\newcommand{\ea}{\end{array}}

\newcommand{\EE}{{\bf E}}
\newcommand{\eps}{\varepsilon}
\newcommand{\sign}{{\rm sign}}

\theoremstyle{plain}
\newtheorem{theorem}{Theorem}[section]
\newtheorem{lemma}[theorem]{Lemma}

\newtheorem{corollary}[theorem]{Corollary}
\newtheorem{remark}[theorem]{Remark}

\newtheorem{example systems}[theorem]{Example systems}

\begin{document}
\title{Validity of the nonlinear Schrödinger approximation for quasilinear dispersive systems with more than one derivative}
\author{Max  He{\ss}\thanks{Universität Stuttgart, Institut für Analysis, Dynamik und Modellierung; Pfaffenwaldring 57, D-70569
Stuttgart, Germany; e-mail: max.hess@mathematik.uni-stuttgart.de}} 
\date{\today}
\maketitle
\begin{abstract}
For nonlinear dispersive systems, the nonlinear Schr\"odinger (NLS) equation can usually be derived as a formal approximation equation describing slow spatial and temporal modulations of the envelope of a spatially and temporally oscillating  underlying carrier wave.
Here, we justify the NLS approximation for a whole class of quasilinear dispersive systems,
{which also includes  toy models for the waterwave problem. 
%with ice cover
}This is the first time that this is done for systems, where a quasilinear quadratic
term is allowed to effectively lose more than one derivative.
{With effective loss we here mean the 
loss still present 
after
  making a diagonalization of the linear part of 
the system such that all linear operators in this diagonalization  have the same regularity properties.}

%We set no bound on the amount of regularity a quadratic quasilinear term is allowed to lose,
%apart from it not being stronger than the linear term of the system.
%For the error estimates, we start with a modified energy based on some normal form transformations, for which occurring resonances also %have to be avoided.
%The energy then gets modified further in order to permit the closing of the error estimates by Gronwall's inequality. 
% We give examples with some double dispersion equations and simple water wave models.
\end{abstract}

%\keywords{NLS approximation, quasilinear dispersive system, double dispersion equation,4}

\section{Introduction}

Nonlinear dispersive systems can  be  very difficult to solve as well analytically as numerically. Thus, a valid NLS approximation can be a great tool for understanding the dynamics of such systems.
An introduction to this theory can be found in \cite{SU17}.
%Here, we justify the NLS  approximation for a whole class of quasilinear dispersive systems. 
For the sake of simplicity, we will restrict ourselves in this introduction to the basic prototype equation
\begin{align}
\label{beam equation}
\partial_t^2 u = -\partial_x^4 u -\partial_x^4 u^2  ,
\end{align}
with  $u:\mathbb{R} \times \mathbb{R} \rightarrow \mathbb{R}: (x,t) \mapsto u(x,t) $,
which  is a simple beam equation. Beside modeling deformations of an elastic beam,
it appears for surface waves in shallow water, in the dislocation theory of crystals or for the interaction between waves guides and some external medium,
cf. e.g. \cite{LG19, KV19, WC06}. 
This equation 
is part of the  class of quasilinear dispersive systems for which we justify the NLS approximation in this paper,
namely
\begin{align}
\label{QDS}
%\tag{QDS}
\partial_t {u} = &   -i \omega {v}\,,
\\ 
\nonumber
%\label{system orig2}\tag{QDS}
\partial_t {v}  = &  -i \omega {u} - i \rho u^2 \,,
\end{align}
with $\rho$ and $\omega$ being differential operators. 
{These systems
are already diagonalized such that the  operators
acting on the linear terms  both cause the same loss of regularity, thus 
 the effective loss caused by the quasilinear terms is the loss caused by $\rho$.}
 For \eqref{beam equation},  the operators $\omega$ and $\rho$ are given in Fourier space by the multipliers $\omega(k)= \rho(k)= \sign(k) k^2$,
 { so \eqref{beam equation} has a quasilinear quadratic term effectively losing two derivatives.}
\\
In order to derive the NLS equation for \eqref{QDS}, we  make an ansatz of the form
\begin{align*}
u&= \varepsilon \psi_{NLS}+ \mathcal{O}(\varepsilon^2) \,,
\end{align*}
where
\begin{align}
\label{psi NLS}
\varepsilon \psi_{NLS}(x,t)= \varepsilon A \big(\varepsilon(x-c_gt), \varepsilon^2 t \big) e^{i(k_0x-  \omega_0 t)}+c.c. \,.
\end{align}
%is the NLS approximation for solutions of \eqref{QDS}.
$0< \varepsilon \ll 1$  is a small perturbation parameter,  $A$ is a complex-valued amplitude and $c.c.$ the complex conjugate.
The ansatz leads to
a basic temporal wave number $\omega_0= \omega(k_0)$ associated to
the basic spatial wave number $k_0>0$ and a group velocity  $c_g= \omega'(k_0)$.
Most importantly, the NLS equation 
\begin{align}
\label{NLS intro}
 \partial_T A = i\, \frac{\omega''(k_0)}{2}\, \partial^2_X A + i\nu_2 (k_0) A |A|^2 ,
\end{align}
is obtained as a lowest order modulation equation,  describing slow modulations in time and space of the envelope of the 
%temporally and spatially oscillating
 wave packet. 
$T= \varepsilon^2 t$ is the slow time scale,  $X= \varepsilon(x-c_gt)$ the slow spatial scale
and $\nu_2 (k_0) \in \mathbb{R}$. A formula for $\nu_2 (k_0)$ is given later.
The NLS equation \eqref{NLS intro} can be explicitly solved, see e.g. \cite{AS81}.

\begin{figure}[H]
\centering
\includegraphics[height=0.5\textwidth ,width=0.9\textwidth]{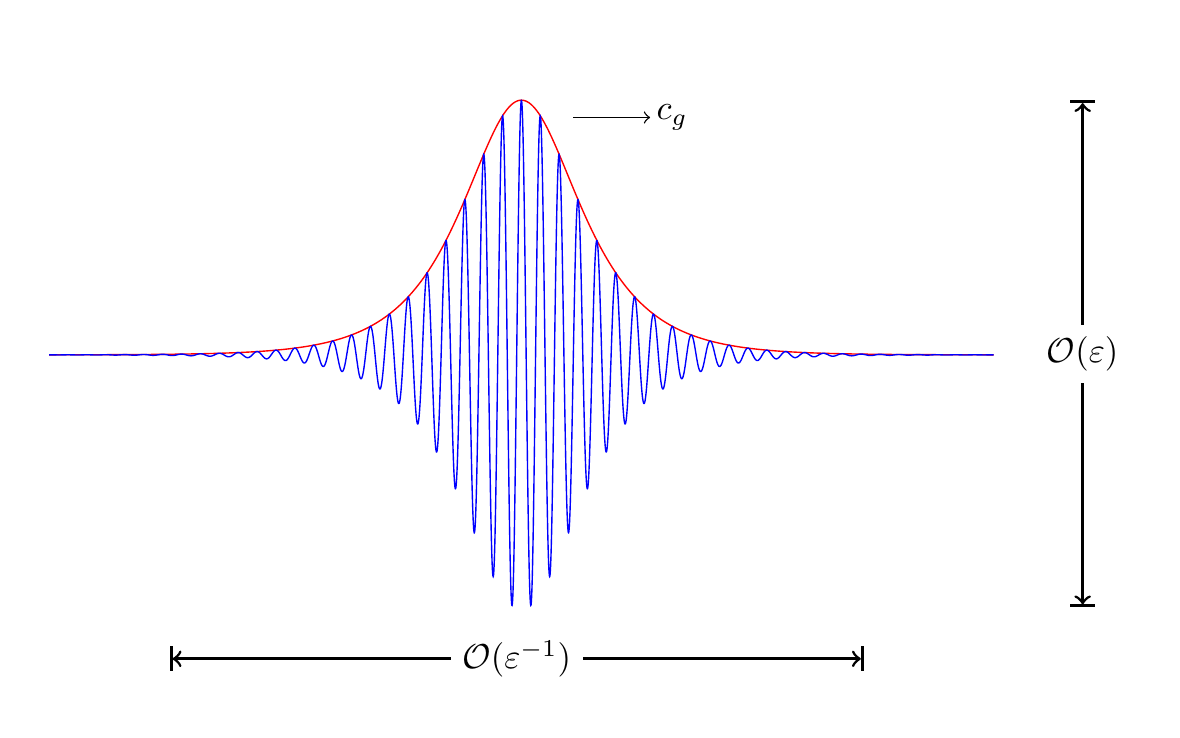}
\caption{The NLS approximation {$\psi_{NLS}$}, an oscillating wave packet with an {envelope}  determined by the solution $A$ of the NLS equation \eqref{NLS intro}. }
\end{figure}
We prove the following result.
\begin{theorem}
\label{mainresult}
Fix $\omega$, $\rho$ and $k_0>0$. 
%\\
For $s_A \geq 7$ and all $C_1,T_0 > 0$ 
there exists some $\varepsilon_0, C_2>0$ 
such that for all solutions $A \in
C([0,T_0],H^{s_A}(\R,\C))$ of the NLS equation \eqref{NLS intro}
with 
\begin{align*}
\sup_{T \in [0,T_0]} \| A(\cdot,T) \|_{H^{s_A}(\R,\C)} \leq C_1
\end{align*}
the following holds.
\\
For all $\varepsilon \in\, (0,\varepsilon_0)$
there are solutions 
\begin{align*}
u \in
C \big([0,T_0/\varepsilon^2], H^{s_A}(\R,\R) \big)
\end{align*}
of the original system \eqref{QDS} which satisfy
\begin{align*}
\sup_{t\in[0,T_0/\varepsilon^2]} \| u(\cdot,t) -
    \varepsilon \psi_{NLS}(\cdot,t)\|_{H^{s_A}(\R,\R)} 
\le  \varepsilon^{3/2} C_2.
\end{align*}
\end{theorem}
Our estimate states that the dynamics of the NLS equation are present in the original system.
Such a result is nontrivial and should never be taken for granted,  see e.g. \cite{SSZ15}.
\\
The following properties of  system \eqref{QDS} make a justification of the NLS equation especially difficult:
\begin{itemize}
\item
a quasilinear quadratic term  in the presence of  nontrivial resonances,
\item
a nonlinearity that causes the loss of more than one derivative in the error estimates.
\end{itemize}
This is the first paper where an NLS justification result is given for a quasilinear system where the nonlinearity { causes an effective loss of} more than one derivative. In this case  a qualitatively new analysis is needed
since a loss of regularity in the error evolution can no longer be dodged through integration by parts like in previous articles.
\\
A quadratic term yields in the equations of the error to terms of order $\mathcal{O}(\varepsilon)$, which could potentially lead to an explosion on a time scale of order  $\mathcal{O}(\varepsilon^{-1})$.
In numerous articles, a theory to handle quadratic terms by using normal form transformations was developed, see e.g. \cite{K88}  and 
for the case of resonances \cite{S05}. However, quasilinear quadratic terms were explicitly excluded. Such terms
make the closing of error estimates much harder.
%On the qualitatively correct timescale, 
The  first NLS validity results  for systems with a quasilinear  quadratic term
(losing a half derivative)
  then were proven in \cite{SW11}
%with the help of a a Cauchy-Kowalevskaya argument only working for quadratic terms that lose not more than half a derivative. In and
and  \cite{TW12}.
%a quasilinear quadratic term losing half a derivative was handled via some special coordinate transform. 
%In \cite{CS13}, numerical evidence for the  validity of the NLS
%in case of a quasilinear term losing a  whole derivative was given.
%{In  \cite{HITW15},  Hunter, Ifrim, Tataru and Wong  developed the idea of normal form
%transforms further by introducing a modified energy which can circumvent the non-invertibility of normal form transformations.}
Using a modified energy to handle the occurring loss of regularity in the normal form %transforms as in \cite{HITW15}, 
%the result of \cite{SW11} was improved in \cite{CW17} and 
{NLS validity results for systems with}
quadratic terms that {effectively} lose a whole derivative could be {proven} 
in \cite{D17,DH18}. 
We here now for the first time
 prove a NLS validity result  for {a class of} systems with quasilinear quadratic terms that can {effectively} lose an arbitrary amount of derivatives.

{It has to be mentioned that the problem of proving NLS justification results
%(in the presence of quasilinear terms) 
has always been closely related to the water wave problem. 
Indeed the first time the NLS equation was derived as a model equation, was for the water wave problem in \cite{Z68}.
The (2-D) water wave problem (WWP) is the problem of finding the irrational flow of an incompressible fluid in an infinitely long canal with flat bottom and a free surface under the influence of gravity. 
For the WWP the NLS equation was  rigorously justified on the right time scale by Totz and Wu for the case of zero surface tension and infinite depth
in \cite{TW12}, and  by  D{\"u}ll, Schneider and Wayne for the case of zero surface tension and finite depth in \cite{DSW16}. In both cases quasilinear terms lose effectively a half derivative. In \cite{IT19}
% the NLS equation  was again justified for
%the WWP in the case of zero surface tension and infinite depth 
a result similar to the one of \cite{TW12} was proven
by using a modified energy method. The  NLS equation for the WWP in case of finite depth and possibly of surface tension
was justified,  uniformly with respect to the strength of the
surface tension as the height of the wave packet and the surface tension go to zero, in 
\cite{D19}. In \cite{D19} it was also shown that  due to cancellations  there  only is an effective loss of one derivative, 
%in the WWP (in the arc-length formulation)
 although one might would expect a loss of one and a half derivatives. A claim that such cancellations always occur for physical systems of relevance would be quite optimistic. Our article now also hopes to contribute to a justification of the WWP in the case of finite depth and an ice cover, in which one expects an effective loss of up to two derivatives. We give more details on how system \eqref{QDS} can be viewed as a toy model for the WWP in the discussions section at the end of the article.
 \\
%Due to the generality of the result, we expect that our  approach will also be useful for more complicated quasilinear dispersive %systems like the  water wave problem with ice cover, where quasilinear terms losing more than one derivative occur.
The plan of the paper is as follows.
\\
{We first quickly explain how the NLS approximation is derived and residual estimates are proven for \eqref{QDS}. Then we justify the NLS approximation by proving error estimates.}
In order to obtain the natural $\mathcal{O}(\varepsilon^{-2})$-time scale of the NLS equation {for the error}, we use 
{a modified energy method, i.e. we use}
a modified energy based on normal form transformations that is equivalent to the squared Sobolev norm of the {error.}
%, similar as in \cite{DH18}.  
The evolution of this energy then {still} contains terms that cannot be estimated without a loss of regularity {since an effective loss of more than one derivative is allowed.} 
{In order to control these terms, we 
 recursively construct
an expression of order $\mathcal{O}(\varepsilon)$,
whose time derivative cancels out with the problematic terms. The construction mainly exploits the
time space relation given for the error and the strict concentration of the NLS approximation in Fourier space around certain multiples of the wavenumber $k_0$.
By adding the constructed $\mathcal{O}(\varepsilon)$-term to our energy,
we can then close the energy estimates such that Theorem \ref{mainresult} follows with Gronwall's inequality.}

{\bf Notation}. 
The 
Fourier transform of a function $u \in L^2(\R,\K)$, with $\K=\R$ or $\K=\C$ is denoted by
$\mathcal{F}(u)(k) = \widehat{u}(k) = \frac{1}{2\pi} \int_{\R} u(x) e^{-ikx} dx$.
$H^{s}(\R,\K)$ is
the space of functions mapping from $\R$ into $\K$,
for which
the norm
$ \| u \|_{H^{s}(\R,\K)} = (\int_{\R} |\widehat{u}(k)|^2 (1+|k|^2)^{s} 
dk )^{1/2} $ 
is finite. The space
$ L^1(s)(\R,\K) $ is defined by $ u \in L^1(s)(\R,\K) \Leftrightarrow u \sigma^s \in L^1(\R,\K)$, where
$ \sigma(x) = (1+x^2)^{1/2}$.
We use $ \lceil \alpha \rceil := \min \{z \in \mathbb{Z}:  z \ge \alpha \}$.
We  write
$
I \le \mathcal{O}(E) \, 
$ for expressions $I$ and $E$,
when there exists some constant $C>0$ such that
$
I \le C \,E\,$.
This constant can then always be chosen independently of $E$ and the small perturbation parameter $\varepsilon$.
%\\
%We sometimes write 
%$
%I = \mathcal{O}(E)  ,
%$
%when we want to express that
%$
%I \le \mathcal{O}(E)$
%%  \qquad \text{
%and
%%} \qquad  
%$- I \le \mathcal{O}(E) .
%$
%Sometimes  we write 
%$I \lesssim E$ instead of  $I \le \mathcal{O}(E) $.

\section{The general class of systems}
\label{allgemeinsystem}
The class of systems for which we  consider the NLS approximation consists of the quasilinear dispersive first order systems \eqref{QDS}
%\begin{align}
%\label{QDS}
%%\tag{QDS}
%\partial_t {u} = &   -i \omega {v}\,,
%\\ 
%\nonumber
%%\label{system orig2}\tag{QDS}
%\partial_t {v}  = &  -i \omega {u} - i \rho u^2 \,
%\end{align}
%with $u:\mathbb{R} \times \mathbb{R} \rightarrow \mathbb{R}: (x,t) \mapsto u(x,t) $
%and $v: \mathbb{R}^2 \rightarrow \mathbb{R}:(x,t) \mapsto v(x,t) $, 
where the pseudo differential operators $\omega$ and $\rho$ can be expressed through some odd real-valued functions $\rho$ and $\omega$ in Fourier space.
%I.e. in Fourier space, we have
% \begin{align*}
%\partial_t \widehat{u}(k,t) & =    -i \omega(k) \widehat{v}(k,t)\,,\\
%\partial_t \widehat{v}(k,t) & =   -i \omega(k) \widehat{u}(k,t) - i \rho(k) (\widehat{u} \ast \widehat{u})(k,t) \,.
%\end{align*}
Such a first order system is also equivalent to a  quasilinear dispersive equation
\begin{equation}
\label{QDE}
%\tag{QDE}
\partial_t^2 {u}= -\omega^2 u- \rho \omega u^2 \,.
\end{equation}
We do not allow the quadratic term of the system to contain more derivatives than the linear one. We express this by  demanding
\begin{align}
\label{rho BED2}
\deg^*(\rho) \le \ \deg(\omega) \,,
\end{align}
where we write $ \deg^*(\gamma) \le  s$ for a function $\gamma:\mathbb{R} \rightarrow \mathbb{R}$ when there exists some constant $C$ such that
$
|\gamma(k)|\le C (1+ |k|)^s  $ for large  $|k|$,
and  $ \deg(\gamma) =  s$ when there also is  some $c>0$ such that 
$
c \,(1+ |k|)^s \le |\gamma(k)|\le C (1+ |k|)^s  $ for large  $|k|$.
One of the functions $\omega$ or $\sign(\cdot) \omega(\cdot)$ as well as
$\rho$ or $\sign(\cdot) \rho(\cdot)$ has to lie in $C^{m_{\omega}}(\mathbb{R})$ for $m_{\omega}:= \max \{5, \lceil \deg(\omega) \rceil +1\}$. In other words, we allow $\omega$ and $\rho$ to have a jump in $k=0$.
We further demand
 that
 \begin{align}
\label{ord nimmt ab rho}
\deg^*(\rho^{(n)}) &\le \deg^*(\rho^{(n-1)})-1 \,,
\\
\label{ord nimmt ab omega}
\deg(\omega^{(n)}) &= \deg(\omega^{(n-1)})-1 \,,
\end{align}
for $n=1,\dots, m_{\omega}$ as long as $\rho^{(n)} \neq 0$, respectively  $\omega^{(n)} \neq 0$.
A  behavior,    typical  for most differential operators.
To guarantee the derivation of the NLS equation \eqref{NLS intro}, we require
\begin{align}
\label{omega''(k0) neq 0}
& \qquad\qquad \qquad\qquad \qquad\qquad \qquad\qquad\omega''(k_0) \neq 0,
\\[2mm]
\label{Bedfuerherleit}
&\omega'(k_0)  \neq \pm \omega'(0)
\quad
\text{and} 
\quad 
 \rho(0)=0,
\qquad \quad 
\text{or}
\qquad \quad 
\lim_{k \rightarrow  0^{+}} \omega(k)=\omega(0^+) \neq 0 \, ,
\\
\label{NR1}
&
\qquad\qquad \qquad\qquad \quad
m \omega(k_0) \neq \pm \omega(m k_0) \, \,  \text{ for } \, m= \pm 2, ... \,, \pm 5 \, 
\end{align}
being fulfilled for the wavenumber $k_0 >0$. 

Real-valued solutions to the equations  
\begin{align}
\label{3 resonanz}
\omega(k)- j_1j_2\omega(k \mp k_0) +j_1 \omega( \pm k_0)=0 \,
\end{align}
with $j_1, j_2 \in \{\pm1\}$, are called resonances.
Resonances are problematic in the presence of quadratic terms and have to be avoided for a well-defined nomal-form transformation. A resonance in $k=k_1$ is called trivial, if the quadratic term vanishes 
for the wavenumber $k_1$ in Fourier space. 
We here restrict us to the case that only a trivial resonance in $k=0$ and nontrivial resonances in $k=\pm k_0$ can occur.
%, i.e. we demand that the only possible solutions to \eqref{3 resonanz} are $k=0$ and $k=\pm k_0$.
Allowing more resonances would be possible by making some adjustments as in \cite{DS06, D19}. 
We naturally also forbid  the left hand side of \eqref{3 resonanz}  converging to zero for $|k| \rightarrow \infty$.
If $ \omega(0^+) \neq  0$, we  additionally demand that 
\begin{align}
\label{Bedingung fuer NFT und j_2=1}
\omega(0^{+}) \neq \pm  2 \omega(k_0)  
\end{align}
or
\begin{align}
\label{BedingungT}
\omega'(k_0) & \neq \pm \omega'(0)
, \,
  \rho(0)=0
\quad
\text {and }  \omega(0^{\pm}) \neq 2 \omega(k_0)+ j \omega(2 k_0) &&  \text{for }j \in \{ \pm 1\} \,
\end{align}
to prevent the happening of further resonances. 

For this class of systems we  prove Theorem \ref{mainresult}.
For \eqref{beam equation}, the conditions \eqref{rho BED2}-\eqref{ord nimmt ab omega} and \eqref{omega''(k0) neq 0}-\eqref{NR1} are obviously true for all $k_0 >0$.
To verify that resonances can only occur  in $k=0$ and $k=\pm k_0$  one makes a case analysis of \eqref{3 resonanz} with the quadratic formula.

\section{Residual estimates}
\label{sec:appr}
All coming calculations get much easier by working with diagonalized  system
\begin{align} \label{diag1}
\partial_t {u}_{-1} (x,t) & =    -i \omega {u}_{-1}(x,t) - \frac{1}{2} i \rho \big({u}_{-1}+{u}_1\big)^{2}(x,t)\,,\\ %\label{diag2} 
\nonumber
\partial_t {u}_1(x,t) & =  i \omega {u}_1(x,t) + \frac{1}{2} i \rho \big({u}_{-1} + {u}_1\big)^{ 2}(x,t) \,,
\end{align}
with ${u}_{-1}(x,t), {u}_{1}(x,t) \in \mathbb{R}$.
One obtains this system from \eqref{QDS}
via the invertible transformation
\begin{equation}  
 \left(
\begin{array}{c} {u}_{-1} \\ {u}_{1}
\end{array}
\right) = 
\frac{1}{2}\left(
\begin{array}{cc}
  1 &  1   \\
 1 & -1    
\end{array}
\right)
\left(
\begin{array}{c} {u} \\ {v}
\end{array}
\right) \, .
 \end{equation}
With an simple ansatz of the form 
\begin{align*}
\left(\begin{array}{c} 
{u}_{-1}
\\
{u}_{-1}\\
\end{array}\right)&=\left(\begin{array}{c} 
 \varepsilon \psi_{NLS}
\\
0 \\
\end{array}\right) 
+ \mathcal{O}(\varepsilon^2)
\left(\begin{array}{c} 
1
\\
1 \\
\end{array}\right) \,,
\end{align*}
one can derive the NLS equation \eqref{NLS intro} by expanding the operators $\omega, \rho$  in Fourier space around integer multiples of the basic wave number $ k_0$ with Taylor's theorem, cf. \cite{DH18}. 
We obtain the NLS equation 
\eqref{NLS intro} with 
\begin{align*}
%\label{nu_2 fall 1}
\nu_2 (k_0)&=-\rho(k_0) \,\Big(  \frac{ \rho(2 k_0) \,\omega(2k_0)}{4 \big(\omega(k_0)\big)^2- \big(\omega(2k_0) \big)^2}+
\frac{2 \rho'(0) \, \omega'(0)}{\big(\omega'(k_0)\big)^2- \big(\omega'(0) \big)^2} \Big) 
\end{align*} 
when $\omega(0^+)=0$
and with\begin{align*}
%\label{nu_2 fall 2}
\nu_2 (k_0)&=- \rho(k_0) \, \Big( \frac{ \rho(2 k_0) \,\omega(2k_0)}{4 \big(\omega(k_0)\big)^2- \big(\omega(2k_0) \big)^2}  -2 \, \frac{\rho(\,0^{+}) }{\omega(\,0^{+})}          \Big) \, 
\end{align*}
when $\omega(0^+)\neq 0$. The explicit computation done here, can be found in \cite{H19}.

The residual $ {\rm Res}_u(\eps \Psi)$ of an approximation $\eps \Psi$ denotes all terms that remain after 
plugging in an approximation  $\eps \Psi$ into the equations of system \eqref{QDS}. For the coming error estimates,  a very small residual whose norm  can be controlled in high Sobolev spaces is needed.
For this reason an improved approximation $\eps\Psi$ is used.
By exploiting that an ansatz like above is always  strongly concentrated around a finite number of integer multiples of the basic wave number $ k_0 > 0 $,
cut-off functions can be used to restrict the support of an ansatz in Fourier space to small neighborhoods of these wave numbers $j k_0$ with $j \in \{-5,\dots,5\}$. This way, an  approximation $\eps\Psi$ that is an analytic function and has a residual of the formal order $\mathcal{O}(\varepsilon^6 )$ is obtained, cf. Section 2 of \cite{DSW16}.
\\
The approximation that we use is
\begin{equation}
\label{ansatz2}
\eps\Psi = \eps\Psi_c +  \eps^2\Psi_q\,,
\end{equation}
where
\begin{eqnarray*}
\eps \Psi_{c} &=& \eps \psi_{c} \, \!\left( \ba{c} 
1 \\ 0  \ea \right) \;\; 
= \;\; \eps (\psi_{1}+ \psi_{-1}) \, \!\left( \ba{c} 
1 \\ 0  \ea \right) \;\; 
\\
&=&\;\;  \eps \, \big( \,  A_1 (\eps
(x-c_g t),\eps^2 t) \,\EE +c.c. \big) \, \!\left( \ba{c} 
1 \\ 0  \ea \right),
\\[3mm]
\eps^2\Psi_q &=& \eps^2  \, \left(\begin{array}{c}
\psi_{q_{-1}} \\
 \psi_{q_{1}}
\\ \end{array}\right)\;\; =\;\;
 \eps^2\Psi_0
 + \eps^2\Psi_2  +  \eps^2\Psi_h\,,
 %\\[3mm]
\end{eqnarray*}
\begin{eqnarray*}
\eps^2 \Psi_0 &=& \;\; \eps^2 \, \left(\begin{array}{c} A_{0} (\eps
(x-c_gt),\eps^2t)\\
D_{0} (\eps (x-c_gt),\eps^2t)\\
\end{array}\right),
\\[3mm]
\eps^2 \Psi_{2} &=& \;\; \eps^2 \, \left(\begin{array}{c} 
 A_{2} (\eps (x-c_gt),\eps^2t)\,\EE^{2} + c.c.\\
 D_2 (\eps
(x-c_gt),\eps^2t)\,\EE^{2} + c.c.\\
\end{array} \right),
%\\[3mm]
\end{eqnarray*}
\begin{eqnarray*}
\eps^2\Psi_h & = &
\sum_{n=1,2,3,4} \eps^{1+n} \, \left(
\begin{array}{c} A^n_{1} (\eps(x -c_gt),\eps^2t)\EE+ c.c.\\
D^{n}_{1} (\eps(x -c_gt),\eps^2t)\EE + c.c.\\
\end{array}\right)\\
&& + \sum\limits_{n=1,2,3} \eps^{2+n} \,
\left(\begin{array}{c}A^n_{0} (\eps(x -c_gt),\eps^2t)\\
D^n_{0} (\eps(x -c_gt),\eps^2t)\\
\end{array}\right)\\
&& + \sum\limits_{n=1,2,3} \eps^{2+n} \,
\left(\begin{array}{c}A^n_{2} (\eps(x -c_gt),\eps^2t)\EE^2 +c.c.\\
D^n_2 (\eps(x -c_gt),\eps^2t)\EE^2 +c.c.\\
\end{array}\right)\\
&& + \sum\limits_{ n=0,1,2}  \eps^{3+n} \,
\left(\begin{array}{c}A^n_{3} (\eps(x -c_gt),\eps^2t)\EE^3 +c.c.\\
D^n_{3} (\eps(x -c_gt),\eps^2t)\EE^3 +c.c.\\
\end{array}\right)\\
&& + \sum\limits_{ n=0,1} \eps^{4+n}
\left(\begin{array}{c} A^n_{4} (\eps(x -c_gt),\eps^2t)\EE^4 +c.c.\\
 D^n_{4} (\eps(x -c_gt),\eps^2t)\EE^4 +c.c.\\
\end{array}\right),
\\
&& + \,  \eps^{5}
\left(\begin{array}{c} A^0_{5} (\eps(x -c_gt),\eps^2t)\EE^5 +c.c.\\
 D^0_{5} (\eps(x -c_gt),\eps^2t)\EE^5 +c.c.\\
\end{array}\right),
\end{eqnarray*}
where $\EE = e^{i(k_0x-  \omega_0 t)} $, $ \omega_0 = \omega(k_0) $ and $c_g =  \omega' \,(k_0) $.
Here, $  A_1 \big(\eps
( \cdot-c_g t),\eps^2 t\big) $ is the restriction of $ A \big(\eps
( \cdot-c_g t),\eps^2 t\big) $  in Fourier space to the interval $ \{k \in \mathbb{R}: |k| \le \delta < k_0/20 \}$ by some cut-off function, while $A$  is the solution of the NLS-equation \eqref{NLS intro} and $\delta>0$. More precisely
\begin{align}
\label{A_1}
A_1\big(\eps
( \cdot-c_g t),\eps^2 t\big)
%&= P_{0, \delta} \big[ A\big(\eps
%( \cdot-c_g t),\eps^2 t\big)  \big]
%\\
&:=
\mathcal{F}^{-1}
\Big[ \chi_{[-\delta, \delta]}(\cdot) \mathcal{F} \big[ A\big(\eps
( \cdot-c_g t),\eps^2 t\big)  \big](\cdot) \Big] \,,
\end{align}
where  %$\mathcal{F}$ is the Fourier transform and 
$\chi_{[-\delta, \delta]}$  is the characteristic function on the interval $[-\delta, \delta]$, i.e.
$\chi_{[-\delta, \delta]}(k)=1$ for $[-\delta, \delta]$ and $\chi_{[-\delta, \delta]}(k)=0$ for $ k \notin [-\delta, \delta]$.
One can think of $\eps\psi_c$ as $\eps\psi_{NLS}$ with a support in Fourier space restricted to small neighborhoods of the wave numbers $\pm k_0$.
The $ A^n_{j}$ and $ D^n_{j}$ are chosen suitably depending on $A_1$ such that the supports of $ A^n_{j}\EE^j$ and $ D^n_{j} \EE^j$ in Fourier space
lie in  small neighborhoods of the wave number $j k_0$.
%\\
%
%
%\begin{remark}
%It is important to mention that we also have $\widehat{\psi}_c (k)=0$ for $|k| \le \delta$, what can be most helpful when dealing with resonances.
%\end{remark}
\\
Similarly as in  \cite{DSW16}, one obtains:
\begin{lemma} \label{lemmaRES}
Let $ s_A \geq {7} $ and
$ A \in C([0,T_0], H^{s_A}(\R,\C)) $ be a solution
of the NLS equation \eqref{NLS intro}  with 
$$ \sup_{T \in [0,T_0]} \| 
A \|_{H^{s_A}} \leq C_A . $$ 
Then for all $s \geq 0$ there exist $ C_{Res}, C_{\Psi}, \eps_0>0  $ depending on  $C_A$ such that for all
$ \eps \in (0,\eps_0) $ the approximation $\eps \Psi= \eps\Psi_c +  \eps^2\Psi_q$ satisfies
\begin{eqnarray} \label{RES1}
\sup_{t \in [0,T_0/\eps^2]} \|  {\rm Res}_u(\eps \Psi)
\|_{H^s}
& \leq & C_{\rm Res}\, \eps^{{11/2}}, \\ 
\label{RES2}
\sup_{t \in [0,T_0/\eps^2]} 
\big\Vert\eps \Psi - \varepsilon \psi_{NLS} 
\left( \begin{array}{c} 1
\\
0 
\end{array} \right)
\big\Vert_{{H^{s_A}}}
& \leq & {C_{\Psi}}\, \eps^{3/2},
\\
\label{RES3}
\sup_{t \in [0,T_0/\eps^2]} (\|\widehat{\Psi}_{c} \|_{L^1({s+1})(\R,\C)}
+ \|\widehat{\Psi}_q \|_{L^1({s+1})(\R,\C)})
&{\leq} & {C_{\Psi}}\,,
\\
 \label{dtpsi}
\| \partial_t  \widehat{\psi}_{\pm 1} + i \widehat{\omega \psi}_{\pm1} \|_{L^{1}(s)} &\le &C_{\Psi}\, \eps^2\,.
\end{eqnarray}
\end{lemma}
Due to the bound \eqref{RES2}, we can work with the improved approximation $\varepsilon \Psi$ to obtain a validity result for the NLS approximation $ \varepsilon \psi_{NLS}$.
\\
The bound \eqref{RES3} will be needed to  make estimates like
\begin{align*}
\| \psi_{c} f \|_{H^{s}} \leq C \| \psi_{c}  \|_{C^{s}_b}\| f 
\|_{H^{s}} \leq  C \| \widehat{\psi}_{c}  \|_{L^1(s)}\| f
\|_{H^{s}} ,
\end{align*}
without losing powers in $ \varepsilon $ as one would with $\| {\psi}_{c}  \|_{H^s}=\| \widehat{\psi}_{c}  \|_{L^2(s)}$,
where the slow space scale of the NLS causes   $\Vert A(\varepsilon \cdot) \Vert_{L^2}= \varepsilon^{-1/2}\Vert A( \cdot) \Vert_{L^2}$.
The bound \eqref{dtpsi} will be used to approximate $\partial_t \psi_{\pm 1}$ for the normal form transforms.

\section{The error estimates}
\label{sec:err}
We 
write the error of the approximation $\varepsilon \Psi$ as
\begin{align}
\label{R}
\varepsilon^{\beta} \left(\!
\begin{array}{c} {\vartheta R_{-1}} \\ \vartheta {R}_1
\end{array} \!
\right)
:=
\left( \!
\begin{array}{c} {u_{-1}} \\ {u}_1
\end{array} \!
\right) -\varepsilon \Psi  ,
\end{align}	
where $R_{-1}$ and $R_{1}$ are our error functions, $\beta = 5/2$ and
$\vartheta$ is an invertible operator given in Fourier space  either by the weight function
\begin{align}
\label{varthetadef}
\hat{\vartheta}(k)= 
\begin{cases}
  \varepsilon +(1 - \varepsilon ) \frac{|k|}{\delta}  & \quad \text{for }  |k| \le \delta \,,\\
    1 & \quad \text{for }  |k| > \delta \,,
\end{cases}
\end{align}
or by $\hat{\vartheta}(k)=1$ if $0 \neq \pm  \omega(0^+) \neq 2 \omega(k_0)$.
The fixed parameter $\delta$ is as in \eqref{A_1}.
\\
The inclusion of the operator $\vartheta$ is essential for our handling of the nontrivial resonances in $k= \pm k_0$.
For this purpose $\vartheta$ has also already been used e.g. in \cite{DS06, DH18}.
When $0 \neq \pm  \omega(0^+) \neq 2 \omega(k_0)$, there is no resonance in $k=k_0$ such that setting $\hat{\vartheta}(k)=1$ is better.
\\
By  plugging in the above definition into the diagonalized system 
we obtain the following
dynamics for the error
\begin{equation}
\label{R_j}
\begin{aligned} 
\partial_t {R}_{-1}  
 &= 
   -i \omega{R}_{-1}
 - \veps i \rho \vartheta^{-1} (R_{\psi}  (\vartheta{R}_{-1}+\vartheta{R}_1))
+ \veps^{-\beta} \vartheta^{-1}  {\rm Res}_{u_{-1}}(\veps \Psi)
\,, 
\\
\partial_t {R}_{1}  
 &= 
   i \omega{R}_{1}
 + \veps i \rho \vartheta^{-1} (R_{\psi}  (\vartheta{R}_{-1}+\vartheta{R}_1))
+ \veps^{-\beta} \vartheta^{-1}  {\rm Res}_{u_{1}}(\veps \Psi) 
\,, 
\end{aligned}
\end{equation}

where 
\begin{align}
\label{Rpsi}
R_{\psi}:=\psi+\frac{1}{2} \varepsilon^{\beta-1} (\vartheta R_{-1}+\vartheta R_{1}),
\\
\label{kleinpsi}
\psi:=\psi_c+ \varepsilon \psi_Q :=\psi_c+ \varepsilon(\psi_{q_{-1}}+\psi_{q_{1}}) \,.
\end{align}
\\

Following the  idea of \cite{DH18}, we define the modified energy 
\begin{align} 
\label{energy}
\mathcal{E}_{\ell} =  E_{0} +E_{\ell}\,, 
\end{align} 
\begin{align*}
E_{\ell} = \sum_{j_1\in \{\pm1\}} \Big(\, \frac{1}{2} \,  \big\Vert \partial_x^{\ell} R_{j_1} \big\Vert_{L^2}^2  +     \veps \sum_{j_2 \in \{\pm1\}} \int_{\R}\partial_x^{\ell} R_{j_1} \partial_x^{\ell} \vartheta^{-1}N_{j_1 j_2}(\psi_c,R_{j_2})\,dx \Big),
\end{align*}
\begin{align*}
%\label{Def_E0}
E_0(R)= \Vert \check{R}_{-1}\Vert_{L^2}^2 +\Vert \check{R}_1\Vert_{L^2}^2 \,, 
\end{align*}
where
\begin{equation*}
\check{R}_j=R_j + \veps \sum_{j_2 \in \{\pm1\}} \vartheta^{-1}N_{j j_2}(\psi_c,R_{j_2}) 
+ \varepsilon^2  \sum_{j_2,j_3 , j_4 \in \{\pm1\}} \vartheta^{-1}\mathcal{T}_{j j_2 j_3 j_4}(\psi_{j_4},\psi_{j_4},R_{j_3}) \,,
\end{equation*}
\begin{align*}
%\label{Def_Nj1j2}
\widehat{N}_{j_1 j_2}(\psi_c,R_{j_2})(k) =  \int_{\R} {n}_{j_1 j_2}(k,k-m,m) \widehat{\psi}_c(k-m) \widehat{R}_{j_2}(m)\,dm\,,  
\end{align*}
\begin{align*}
%\label{Def_Tj_1j2j3j4j}
&\widehat{\mathcal{T}}_{j_1 j_2 j_3 j_4}(\psi_{j_4},\psi_{j_4},R_{j_3})(k)= 
\int_{\R} {t}_{j_1, j_2, j_3 , j_4} (k) \widehat{\psi}_{j_4}(k-m) \widehat{\psi}_{j_4}(m-n) \widehat{R}_{j_3}(n)\, dn \,dm\, ,
\end{align*}
\begin{align*}
{n}_{j_1 j_2}(k,k-m,m) = \frac{ \rho(k)  
\,  \hat{\vartheta}_{\varepsilon, \infty}(m) \, {\chi}_c(k-m)}{\omega(k)-j_1 j_2\omega(m) +j_1\omega(k-m)}\,,
\end{align*}
\begin{align*}
&{t}_{j_1, j_2, j_3 , j_4}(k) = 
\dfrac{ -j_2 \,\hat{P}_{0, \delta}(k) \,  n_{j_1 j_2}(k,j_4k_0,k- j_4 k_0)    \, \rho(k- j_4 k_0) }{\big(-j_1\omega(k)-2\omega( j_4 k_0)+j_3\omega(k -2 j_4 k_0)\big)} \,,
\end{align*}
\begin{align*}
\hat{\vartheta}_{\varepsilon, \infty}(m)= 
%\begin{cases}
% \qquad 1 
%&,  \text{when }  0 \neq \pm  \omega(0^+) \neq 2 \omega(k_0),
%\\[4mm]
% \quad
\begin{cases}
0 & \quad \text{for }  |m| \le \varepsilon \, ,\\
 \varepsilon +(1 - \varepsilon ) \frac{|m|}{\delta}  & \quad \text{for }  \varepsilon < |m| \le \delta \, ,\\
    1 & \quad \text{for }  |m| > \delta \, .
\end{cases} 
\end{align*}
We set
${t}_{j_1, j_2, j_3 , j_4}=0$ and $\hat{\vartheta}_{\varepsilon, \infty}= 1$
when  $0 \neq \pm  \omega(0^+) \neq 2 \omega(k_0)$.

%We  now  first 
%confirm the required time scale and the equivalence of the energy to the squared Sobolev norm of %the error.
%
%
%show that evolution of this energy has the $\varepsilon$-order required for the NLS.
%In section \ref{SEC Energy equivalence} we then verify the equivalence of this energy to the squared Sobolev norm of the error.
The final handling  of the loss of regularity stemming from quasilinear quadratic terms losing more than one derivative
is done in section \ref{SEC Closing the energy estimates},  by transforming this energy.
After the transformation, the energy estimates close and the theorem is proven.

\subsection{Natural time scale of the NLS}
In order to achieve error estimates valid on the natural $\mathcal{O}(\varepsilon^{-2})$-time scale of the NLS equation,
the evolution of the energy has to be of order $\mathcal{O}(\varepsilon^2)$.
To met this goal, we chose
 the normal form transformations such that  the $\mathcal{O}(\varepsilon)$-terms  in the error evolution are eliminated.

The operator $\vartheta^{-1}$ was placed outside of the normal form transforms ${N}_{j_1 j_2}$ and $\mathcal{T}_{j_1 j_2 j_3 j_4}$
since
\begin{align}
\label{vartheta⁻1=eps⁻1}
\hat{\vartheta}^{-1}(k)= \frac{1}{\hat{\vartheta}(k)} = \mathcal{O}(\varepsilon^{-1})\,.
\end{align}
In the presence of a spatial derivative, better estimates become possible.
\begin{lemma}
We have
\begin{align}
\label{theta^-1 k}
|k  \hat{\vartheta}^{-1}(k) |\le 1+|k| \,.
\end{align}
In particular
\begin{align}
\label{theta-1rho}
\Vert i \rho \vartheta^{-1} f\Vert_{L^2} &\le \mathcal{O}( \Vert  f \Vert_{H^{\deg^* \rho}}).
\end{align}
%In the case  $\vartheta_{\varepsilon, \infty} \neq id_{L^2}$, there is  a constant $C=C(\delta)$ such that for all $k \in \R$:
%\begin{align}
%\label{theta tanh^-1}
%|k^{-1} \, \hat{\vartheta}_{\varepsilon, \infty}(k)| \le  C.
%\end{align}

\end{lemma}
\textbf{Proof.}
The lemma is obviously true for $\hat{\vartheta}(k)=1$. Otherwise, we have
$|k \, \hat{\vartheta}^{-1}(k)|= |k|$ for $|k| > \delta$ and 
%\begin{align*}
%|k \, \hat{\vartheta}^{-1}(k)|= 
%\begin{cases}
%|k | & \quad {\rm for }\; |k| > \delta\,,\\[1mm]
%\dfrac{|k|}{\varepsilon +(1- \varepsilon) \frac{|k|}{\delta}} & \quad{\rm for}\; |k| \le \delta \,.
%\end{cases}
%\end{align*}
%For $0< |k|  \le \delta$, we have
\begin{align*}
|k \, \hat{\vartheta}^{-1}(k)|=
\frac{|k|}{\varepsilon +(1- \varepsilon) \frac{|k|}{\delta}} 
\le \frac{1}{\frac{\varepsilon}{|k|} + \frac{1- \varepsilon}{\delta}} 
\le \delta
\end{align*}
for $0< |k|  \le \delta$ such that \eqref{theta^-1 k} is true.
When $\hat{\vartheta}^{-1} \neq 1$
we are in the case   $\rho(0)=0$. Thus $\rho(k) = \mathcal{O}(k)$ for $|k| \rightarrow 0$ such that \eqref{theta-1rho}  follows.
%When $\vartheta_{\varepsilon, \infty} \neq id_{L^2}$, we have 
%$|k^{-1}\, \hat{\vartheta}_{\varepsilon, \infty}(k)| =|k|^{-1} \le \delta^{-1}$ for $ |k| \ge \delta$,
%\begin{align*}
%|k^{-1}\, \hat{\vartheta}_{\varepsilon, \infty}(k)|
%= 
%%\left\{ \begin{array}{ll}
%% 0  & \quad {\rm for }\; 0 < |k| \le \varepsilon\,, \\[2mm]
%\dfrac{\varepsilon}{|k|} + \dfrac{(1 - \varepsilon )}{\delta } 
%% & \quad {\rm for }\; \varepsilon \le |k| \le \delta\, , \\[4mm]
% % \dfrac{1}{|k|} & \quad {\rm for }\;  |k| \ge \delta \,.
%  %\end{array} \right.
%  \le C
%  \,,
%\end{align*}
%for $\varepsilon \le |k| \le \delta$ and $|k^{-1}\, \hat{\vartheta}_{\varepsilon, \infty}(k)| =0$ for $ |k| \le \varepsilon$. Thus \eqref{theta tanh^-1} is true.
\qed
\medskip

\begin{lemma}
\label{lemma Normalformtransfo}
The normal form transforms $N_{j_1 j_2}$ were constructed such that
for all $f \in H^{ \deg^*(\rho)+1}(\R)$: 
\begin{align} \label{nf}
-j_1 i \omega N_{j_1 j_2}(\psi_c,f)- N_{j_1 j_2}(i \omega \psi_c,f)+ j_2 N_{j_1 j_2}(\psi_c,i \omega f)  &=
-j_1  i \rho ( \psi_c \vartheta_{\varepsilon, \infty} f) \,,
\end{align}
where
\begin{align}
\label{nf-rest}
\varepsilon \,\Vert j_1  i \rho \vartheta^{-1}( \psi \vartheta  f) -j_1  i \rho \vartheta^{-1} ( \psi_c \vartheta_{\varepsilon, \infty} f) \Vert_{L^2}
=
\mathcal{O}(\varepsilon^2) \, \Vert f \Vert_{ H^{ \deg^*(\rho)}}.
\end{align}
%i.e.
%\begin{align*} 
%-j_1 i \omega N_{j_1 j_2}(\psi_c,f)- N_{j_1 j_2}(i \omega \psi_c,f)+ j_2 N_{j_1 j_2}(\psi_c,i \omega f)  +j_1  i \rho ( \psi_c \vartheta  f)&=
%  \mathcal{O}(\varepsilon).
%\end{align*}
Moreover, the operators $N_{j_1 j_2}(h,\cdot)$ are continuous linear operators which map $H^1(\R,\R)$ into $L^2(\R,\R)$
for fixed $h \in L^2(\R,\R)$.
In particular, there is a $C=C(\Vert \widehat{h}(\cdot) \chi_c(\cdot) \Vert_{L^1})$ such that for all $g \in H^1(\mathbb{R)}$:
\begin{align}
\label{Njj absch}
\Vert N_{j j}(h,g) \Vert_{L^2} &\le C \Vert g \Vert_{H^1} \,,
\\[1mm]
\label{Nj-j absch}
\Vert N_{j -j}(h,g) \Vert_{L^2} &\le C  \Vert g \Vert_{L^2}.
\end{align} 
\end{lemma}
{\bf Proof.}
In order prove that the $N_{j_1 j_2}$ are well-defined, we have to look at the zeros of the denominator of ${n}_{j_1 j_2}$, i.e. of
\begin{equation*}
  \omega(k)-j_1 j_2\omega(m) +j_1\omega(k-m)
\end{equation*}
for $|k-m \mp k_0|\le \delta$.
Due to the assumption for \eqref{3 resonanz} in section \ref{allgemeinsystem}, we can chose $\delta$ such small that for $|k-m \mp k_0|\le \delta$ the equation
\begin{align}
\omega(k)-j_1 j_2\omega(m) +j_1\omega(k-m)=0
\end{align}
can have no other solutions than $k=0$ or $m=0$.

We first check $k=0$ and therefore assume $|k| \le \delta$.
\\
For $|k| \le \delta$, we also have $|-m \mp  k_0 |\le 2\delta$ since $|k-m \mp k_0|\le \delta$.
Using Taylor  in order to expand $ \omega(k)$ in the point $ \sign(k)\cdot 0^+$
and  $\omega(k-m)$ in the point $-m$, we obtain
\begin{align*}
&\omega(k)-j_1 j_2\omega(m) +j_1\omega(k-m) 
\\[1mm]
& \qquad\qquad
= \omega(\sign(k) \cdot 0^+)-j_1 j_2\omega(m) +j_1\omega(-m)
\\
&\qquad \qquad \quad
+ \omega'(\sign(k) \cdot 0^+) \, k+ j_1 \omega'(-m) \, k +\mathcal{O}(k^2)
\\[1mm]
&\qquad \qquad
= \omega(\sign(k) \cdot 0^+)-j_1( j_2+1) \, \omega(m) + \big(\omega'(0)+ j_1 \omega'(m) \big) \, k +\mathcal{O}(k^2).
\end{align*}
Thus, if
\begin{align}
\label{BED NFT}
\omega(0^{\mp}) \neq ( j_2+1) \, \omega(k_0),
\end{align}
and we choose $\delta$ small enough,  $N_{j_1 j_2}$ has no resonance in $k=0$.
\\
If \eqref{BED NFT} is hurt but 
\begin{align}
\label{BED2 NFT}
\pm \omega'(0) \neq \omega'(k_0),
\end{align}
we can choose $\delta$ small enough such that
\begin{align*}
\omega(k)-j_1 j_2\omega(m) +j_1\omega(k-m)= \mathcal{O}(k) \qquad \text{for } k \rightarrow 0.
\end{align*}
When \eqref{Bedingung fuer NFT und j_2=1} is true, we have \eqref{BED NFT} and thus
$N_{j_1 j_2}$ has no resonance in $k=0$.
\\
When instead \eqref{BedingungT} is  true, we always have \eqref{BED2 NFT} and $\rho(k)= \mathcal{O}(k)$ for $k \rightarrow 0$
, thus
$N_{j_1 j_2}$  can at worst have a trivial resonance in $k=0$.
%\medskip
%Checking $m=0$ can be done completely analogously since
%\begin{align}
%\omega(k)-j_1 j_2\omega(m) +j_1\omega(k-m)= -j_1 j_2 \big( \omega(m)-j_1 j_2\omega(k) -j_2 \omega(k-m) \big).
%\end{align}
\\
The case $m=0$ works analogously due to symmetry and the choice of $\vartheta$.
To give more details, 
in the problematic case
\begin{align*}
\omega(k)-j_1 j_2\omega(m) +j_1\omega(k-m)= \mathcal{O}(m) \qquad \text{for } m \rightarrow 0 \,
\end{align*}
there occur no nontrivial resonances or a loss of $\varepsilon$-powers, since 
%$|m^{-1}\, \hat{\vartheta}_{\varepsilon, \infty}(m)| =0$ for $ |m| \le \varepsilon$ and
\begin{align*}
|m^{-1}\, \hat{\vartheta}_{\varepsilon, \infty}(m)|
%= \dfrac{\varepsilon}{|m|} + \dfrac{(1 - \varepsilon )}{\delta } 
  \le 1+\delta^{-1}
  \,.
\end{align*}

Resonances for $|k|, |m| \rightarrow \infty$ were excluded in section \ref{allgemeinsystem}.

The property \eqref{nf} can be easily checked in Fourier space. 

Concerning estimate \eqref{nf-rest},
\begin{align*}
\big\Vert j_1  i \rho \vartheta^{-1}( \psi \vartheta  f) -j_1  i \rho \vartheta^{-1}( \psi_c \vartheta_{\varepsilon, \infty} f) \big\Vert_{L^2}
&=
\big\Vert   i \rho \vartheta^{-1}( \varepsilon \psi_Q \vartheta  f)+ i \rho \vartheta^{-1}\big( \psi_c (\vartheta- \vartheta_{\varepsilon, \infty})  f \big) \big\Vert_{L^2}
\\[1mm]
&=
\mathcal{O}(\varepsilon) \, \Vert f \Vert_{ H^{ \deg^*(\rho)}} \,,
\end{align*}
particularly due to \eqref{theta-1rho} and $(\hat{\vartheta}- \hat{\vartheta}_{\varepsilon, \infty})\le \mathcal{O}(\varepsilon)$.

We now will show that the $N_{j_1 j_2}(h,\cdot)$  are continuous linear operators.
\\
For later purposes, we will especially focus on writing the bilinear operators $N_{j_1 j_2}(\cdot,\cdot)$ as a sum of products of linear operators, plus
some smoothing bilinear operator. 
\\
We first look at $N_{jj}$.
\\
For $|k| \rightarrow \infty$, we have
\begin{align*}
n_{jj}(k,k-m,m)= \frac{\rho(k) \, \chi_{c}(k-m)}{ \omega(k) -\omega(m)+j \, \omega(k-m)}.
\end{align*}
We want a form of $n_{jj}(k,k-m,m)$ for $|k| \rightarrow \infty$ that only consists of terms that are products of functions in one variable, plus some smoothing term.
In order to obtain this, we have to examine the denominator.
%\begin{align*}
%\frac{ \chi_{c}(k-m)}{ \omega(k) -\omega(m)+j \, \omega(k-m)}.
%\end{align*}
Using Taylor, we get 
\begin{align*}
\omega(k) -\omega(m)=\omega'(m)\, (k-m)+ r(k,k-m,m),
\end{align*}
where
\begin{align*}
r(k,k-m,m) \chi_{c}(k-m)= \Big( \sum_{l=2}^p \frac{1}{l!} \, \omega^{(l)}(m) \, (k-m)^l + \mathcal{O} \big(\omega^{(p+1)}(m) \big) \, \Big) \chi_{c}(k-m),
\end{align*}
for some sufficiently large chosen $p \ge \lceil \deg^*(\rho) \rceil$.
Then we use the expansion
\begin{align}
\label{bruchentwicklung}
\frac{a}{b+c}= \sum_{l=0}^{n} (-1)^l \frac{a c^l}{b^{l+1}}+ (-1)^{n+1} \frac{ac^{n+1}}{b^{n+1} (b+c)} 
\qquad  \qquad
 (b+c \neq 0, \, b \neq 0).
\end{align}
%in order to obtain a form of 
%\begin{align*}
%\frac{ \chi_{c}(k-m)}{ \omega(k) -\omega(m)+j \, \omega(k-m)} 
%%\quad (\text{for } |k| \rightarrow \infty)
%\end{align*}
%for $|k| \rightarrow \infty$
%which only consists of terms that are products of functions in one variable, plus some
% %$\mathcal{O}(|m|^{-\deg^*(\rho)- \deg^*(\rho')})$-term.
% smoothing term.
%\\
We distinguish the three cases $\deg(\omega)>1$, $\deg(\omega)=1$ and $\deg(\omega)<1$.
\\
If $\deg(\omega)>1$ (i.e. $\deg(\omega') >0$), we have for $|k| \rightarrow \infty$:
\begin{align}
\label{as jj w>1}
&\frac{ \chi_{c}(k-m)}{ \omega(k) -\omega(m)+j \, \omega(k-m)}
\\[2mm] \nonumber
& \qquad
=\frac{ \chi_{c}(k-m)}{\omega'(m)\, (k-m)+ r(k,k-m,m)+j \, \omega(k-m)}
\\[2mm] \nonumber
& \qquad =
 \Big(
\frac{1 }{ \omega'(m)\, (k-m)}
-
\frac{  r(k,k-m,m)+j \, \omega(k-m) }{ \omega'(m)^2\, (k-m)^2}
\\ \nonumber
& \qquad \qquad
+
\frac{\big( r(k,k-m,m)+j \, \omega(k-m) \big)^2}{ \omega'(m)^3\, (k-m)^3}
-
\frac{ \big( r(k,k-m,m)+j \, \omega(k-m) \big)^3}{ \omega'(m)^4\, (k-m)^4}
\\ \nonumber
&  \qquad \qquad
\pm \dots
+
\mathcal{O}(|m|^{-\deg^*(\rho)- \deg^*(\rho')}) \, \Big) \chi_{c}(k-m)
.
\end{align}

If $\deg(\omega)=1$ (i.e. $\deg(\omega') =0$), we have 
\begin{align}
\label{as jj w=1}
&\frac{ \chi_{c}(k-m)}{ \omega(k) -\omega(m)+j \, \omega(k-m)}
\\ \nonumber
& \quad\qquad =
 \Big(
\frac{1 }{ \omega'(m)\, (k-m)+j \omega(k-m)}
+
\mathcal{O}(|m|^{-1)}) \, \Big) \chi_{c}(k-m)
\,,
\quad
\text{for }|k| \rightarrow \infty
.
\end{align}

If $\deg(\omega)<1$ (i.e. $\deg(\omega') <0$), there is some $N=N(\omega') \in \mathbb{N}$ 
such that 
\begin{align}
\label{as jj w<1}
&\frac{ \chi_{c}(k-m)}{ \omega(k) -\omega(m)+j \, \omega(k-m)}
\\ \nonumber
& \qquad =
\Big(
\sum_{n=0}^{N} (-1)^n j^{n+1}
\frac{ (\omega'(m))^n \,(k-m)^{n} }{ ((\omega(k-m))^{n+1}}
+
\mathcal{O}(|m|^{-1}) \, \Big) \chi_{c}(k-m)
\,,
\quad
\text{for }|k| \rightarrow \infty
.
\end{align}
Due to \eqref{rho BED2}, \eqref{ord nimmt ab rho} and \eqref{ord nimmt ab omega}
we now get that
the $N_{j j}(h,\cdot)$ map $H^1(\R)$ on $L^2(\R)$ by taking advantage of
Plancherel's theorem and Young's inequality for convolutions
\begin{align*}
&\Vert N_{j j}(h,g) \Vert_{L^2}
\lesssim  
\Vert \widehat{N}_{j j}(h,g) \Vert_{L^2}
= \big\Vert
\int_{\R} {n}_{j j}(\cdot, \cdot-m,m) \widehat{h}(\cdot-m) \widehat{g}(m)\,dm
\big\Vert_{L^2}
\\[2mm]
&
\quad
 \le
\mathcal{O} \Big( \sup_{k,m \in \mathbb{R}} \! \frac{|{n}_{j j}(k, k-m,m)|}{(|m|^2+1)^{1/2}}  \Big) 
\big\Vert
\int_{\R} | \widehat{h}(\cdot-m)\chi_c(\cdot-m) \, (|m|^2+1)^{1/2} \, \widehat{g}(m)|\,dm
\big\Vert_{L^2}
\\[2mm]
& 
\quad
\le \mathcal{O} \big(\Vert \widehat{h}(\cdot) \chi_c(\cdot) \Vert_{L^1} \big) \, \Vert g \Vert_{H^1} \, .
\end{align*}

Now, we look at $N_{j,-j}$.
\\
Using Taylor, we get for $|k| \rightarrow \infty$:
\begin{align*}
n_{j,-j}(k,k-m,m)
&= \frac{\rho(k) \, \chi_{c}(k-m)}{ \omega(k) +\omega(m)+j \, \omega(k-m)}
\\[3mm]
&
=
\frac{\rho(k) \, \chi_{c}(k-m)}{ 2\omega(k) +r(k,k-m,m) +j \, \omega(k-m)}
,
  \end{align*}
where $r(k,k-m,m)$ is now given by
\begin{align*}
r(k,k-m,m) =  \sum_{l=1}^p \frac{(-1)^{l}}{l!} \, \omega^{(l)}(k) \, (k-m)^l + \mathcal{O} \big(\omega^{(p+1)}(k) \big).
\end{align*}
for some sufficiently large chosen $p \ge \lceil \deg^*(\rho) \rceil$.
Using expansion \eqref{bruchentwicklung},
we obtain 
\begin{align}
\label{asyj-j w al}
n_{j,-j}(k,k-m,m)
&= \Big(
\frac{\rho(k) }{ 2\omega(k) }
-
\frac{\rho(k) \, \big( r(k,k-m,m)+j \, \omega(k-m) \big) }{ 4\omega(k)^2}
\\ \nonumber
& \qquad
+
\frac{\rho(k)\, \big( r(k,k-m,m)+j \, \omega(k-m) \big)^2 }{ 8\omega(k)^3 }
\\ \nonumber
& \qquad
\mp \dots
+
\mathcal{O}(|k|^{-\deg^*(\rho)}) \, \Big) \chi_{c}(k-m) \,,
\qquad \qquad \qquad
\text{for }|k| \rightarrow \infty
.
  \end{align}
We can now see that the $N_{j -j}(h,\cdot)$ map $L^2(\R)$ on $L^2(\R)$ by exploiting Young's inequality for convolutions.

Finally, since
\begin{align*}
n_{j_1 j_2}(-k,-(k-m),-m)=n_{j_1 j_2}(k,k-m,m) \in \mathbb{R},
\end{align*}
the $N_{j_1 j_2}(h,\cdot)$ map real-valued functions on real-valued functions.

\qed

\begin{lemma}
\label{Tlem}
The normal form transforms $\mathcal{T}_{j_1 j_2 j_3 j_4}$ were constructed such that for all $j_1, j_2 , j_3, j_4 \in \{\pm 1\}$, we have 
\begin{align}
\label{T-rest}
\varepsilon^2 \,
\Vert\vartheta^{-1} Y_{j_1, j_2, j_3} \Vert_{L^2} 
\le\varepsilon^2 \, \mathcal{O} \big( \, \Vert R_{j_3} \Vert_{H^{\deg^*(\rho)+1}} \big).
\end{align}
where
\begin{align}
\label{T-wahl}
 Y_{j_1, j_2, j_3}=
&
\,   N_{j_1 j_2 }(\psi_{c},j_2 \vartheta^{-1}\, i \rho (\psi \vartheta R_{j_3}))
  \\[2mm] \nonumber
   &
 +\sum_{j_4= \pm 1} \Big(
     - j_1 \,i \omega  \mathcal{T}_{j_1 j_2 j_3 j_4}(\psi_{j_4},\psi_{j_4},R_{j_3})+
 \mathcal{T}_{j_1 j_2 j_3 j_4}(-i \omega  \psi_{j_4},\psi_{j_4},R_{j_3}) 
 \\[2mm] \nonumber
&
\qquad \qquad \quad
+ \mathcal{T}_{j_1 j_2 j_3 j_4}(\psi_{j_4},-i \omega \psi_{j_4},R_{j_3}) 
+ \mathcal{T}_{j_1 j_2 j_3 j_4}(\psi_{j_4},\psi_{j_4}, j_3 \,i \omega R_{j_3}) \Big) \,.
\end{align}

Furthermore, for fixed functions $g,h$ with $ \widehat{g},\widehat{h} \in L^1(\mathbb{R}, \mathbb{C})$, the mapping
$f \mapsto \mathcal{T}_{j j_3}(g,h,f)$ defines a continuous linear map from $L^2(\R,\C)$ into $L^2(\R,\C)$ and there is a constant 
$C= C \big(\Vert \widehat{g} \Vert_{L^1} \Vert \widehat{h} \Vert_{L^1} \big)$ such that for all $f \in L^2(\R,\C)$, we have
\begin{align}
\label{T-abschaetzung}
\Vert  \mathcal{T}_{j_1 j_2 j_3 j_4}(g,h,f) \Vert_{L^2} \le C  \Vert  f\Vert_{L^2}\,.
\end{align}

%
%{\bf c)}
%For all $f \in L^2(\R,\C)$ we have
%\begin{align}
%\label{P_0 T}
%P_{\delta,\infty}\mathcal{T}_{j_1 j_2 j_3 j_4}(\psi_c, \psi_c,f)=0\,.
%\end{align}
\end{lemma}

%\begin{remark}
%If for  pair $(j_1',j_2')$ the normal form transform $N_{j_1' j_2'}$ has no resonances, 
%all statements of lemma \ref{Tlem} hold true for $\mathcal{T}_{j_1' j_2' j_3}=0$.
%\end{remark}

\textbf{Proof.}
In the case  $0 \neq \pm  \omega(0^+) \neq 2 \omega(k_0) $,  the lemma is trivial.
Otherwise, we have to first look at the zeros of the denominator of $t_{j_1, j_2, j_3 , j_4}(k)$, i.e. the zeros of 
\begin{align*}
\big(\omega(k)- j_1 j_2\omega( k- j_4 k_0)+j_1\omega(j_4 k_0) \big) \, \big(
-j_1\omega(k)-2\omega( j_4 k_0)+j_3\omega(k -2 j_4 k_0)\big)
\end{align*}
for $|k|\le \delta$.
For the first factor, we have \eqref{3 resonanz}, so the only possible zero of the first factor is $k=0$.
For the second factor, we get by expanding the expression $\omega(k)$ in the point $\sign(k) \cdot 0^+$ and $\omega(k -2 j_4 k_0)$ in the point $-2 j_4 k_0$:
\begin{align*}
&-j_1\omega(k)-2\omega( j_4 k_0)+j_3\omega(k -2 j_4 k_0)
\\
&
\qquad
= -j_1\omega(\sign(k) \cdot 0^+)-2\omega( j_4 k_0)+j_3\omega(-2 j_4 k_0)+ \mathcal{O}(k).
\end{align*}
We can choose $\delta$ such small that this expression has no zeros.
For  $\omega(0)=0$ this is  possible due to \eqref{NR1}, and for  $\omega(0) \neq 0$ due to \eqref{BedingungT}. 
Summing up, there can only occur a trivial resonance in $k=0$.
We now  obtain \eqref{T-abschaetzung} by using Young's inequality for convolutions 
and the fact that $ \|\widehat{t}_{j_1, j_2, j_3 , j_4}\|_{L^{\infty}}$ can be uniformly bounded.
%\begin{align*}
%\Vert  \mathcal{T}_{j_1 j_2 j_3 j_4} (g,h,R_{j_3}) \Vert_{L^2} &\le   \|t_{j_1, j_2, j_3 , j_4}\|_{L^{\infty}} 
% \Vert \widehat{g} \Vert_{L^1}
%  \Vert  \widehat{h}  \Vert_{L^1} \Vert  R_{j_3}\Vert_{L^2}
%\le   C  \Vert  R_{j_3}\Vert_{L^2}.
%\end{align*}
\\
The estimate \eqref{T-rest} is obtained similarly as in \cite{DH18}, the  details can be found in \cite{H19}.
\qed
\medskip

\begin{lemma}
\label{lemma d/dt E_0}
For  $m \ge \deg^*(\rho)+1$, we have
\begin{align}
\label{d/dt E_0}
\partial_t E_0 \le 
\eps^2 \, \mathcal{O} \Big(\varepsilon^{1/2} \, \big( \Vert R_{-1} \Vert_{H^m}^2
 + \Vert R_{1} \Vert_{H^m}^2\big )^{3/2}+ \Vert R_{-1} \Vert_{H^m}^2
 + \Vert R_{1} \Vert_{H^m}^2+1 \Big) \,.
\end{align}
\end{lemma}
\textbf{Proof.}
The statement follows by the construction of the normal form transforms,
i.e.  with \eqref{dtpsi} and the lemmata \ref{lemma Normalformtransfo} and \ref{Tlem}. 
%, \eqref{nf} and \eqref{nf-rest}, and,  \eqref{T-rest} and  \eqref{T-wahl}.
A complete proof can be found in \cite{H19}.
\qed

\subsection{Energy equivalence}
\label{SEC Energy equivalence}
In order to obtain a result for the error, the energy $\mathcal{E}_{\ell}$ has to be equivalent to the $H^{\ell}$-energy of the error functions.
{A loss of regularity caused by the normal form transforms
is avoided by using the modified energy $\mathcal{E}_{\ell}$. The first time a modified energy was used to overcome a loss of regularity was in \cite{HITW15}.}
% We have to overcome the loss of regularity caused by the normal form transforms.

\begin{lemma}
\label{lemma E_0}
There are constants $C_0, \check{C}_0$ such that
%, for $\varepsilon \le \varepsilon_0$ and $\varepsilon_0$ small enough, 
the following estimates hold 
\begin{align}
\label{checkR gegen R}
%\Vert \check{R}_1 \Vert_{L^2} +\Vert \check{R}_{-1} \Vert_{L^2} 
\sqrt{E_0}
\le C_0  \,\big( \Vert R_1 \Vert_{H^{1}}+\Vert R_{-1} \Vert_{H^{1}} \big) ,
\end{align}
\begin{align}
\label{R gegen checkR}
\Vert R_1 \Vert_{L^2} +\Vert R_{-1} \Vert_{L^2}\le \check{C}_0 
%\big(\Vert \check{R}_1 \Vert_{L^2} +\Vert \check{R}_{-1} \Vert_{L^2} \big)
\, \sqrt{E_0}
+ \varepsilon \, \mathcal{O} ( \Vert {R_{-1}} \Vert_{L^2} +\Vert {R_1} \Vert_{L^2} ).
%+\mathcal{O}(\varepsilon) \big( \Vert {R_{-1}}\Vert_{H^{1}}+ \Vert {R_1}\Vert_{H^{1}}\big). 
\end{align}
\end{lemma}
\textbf{Proof.}
The proof is similar to the one in \cite{DH18} and can be found in \cite{H19}.
\qed
\medskip

\begin{lemma}
\label{lemma verallgemeinerte partielle Integration}
Let $f,g,h\in L^2(\R,\R)$ be real-valued functions and $K:\R ^3 \rightarrow \C$.
%\textbf{a)}
\\
If
\begin{align*}
\int_{\R} \int_{\R}\big\vert K(k,k-m, m) \;
 \overline{\widehat{f}(k)} \, \widehat{h}(k-m)\, \widehat{g}(m)\big\vert \,dm\,dk < \infty,
\end{align*}
then we have
\begin{align}
\label{veralg PI in FR}
&\int_{\R} \int_{\R} K(k,k-m, m)\;
 \overline{\widehat{f}(k)} \, \widehat{h}(k-m)\, \widehat{g}(m)\,dm\,dk
\\ \nonumber
&  \quad \quad \quad=
 \int_{\R} \int_{\R} K(-m,k-m, -k)\;
 \overline{\widehat{g}(k)} \, \widehat{h}(k-m)\, \widehat{f}(m)\,dm\,dk
 .
\end{align} 
\end{lemma}
\textbf{Proof.}
The result is obtained by exploiting the fact that 
$\overline{\widehat{f}(k)}=\widehat{f}(-k)$ and $\widehat{g}(m)=\overline{\widehat{g}(-m)}$,
making a change of variables  and  using Fubini's theorem.
%\begin{align*}
%&\int_{\R} \int_{\R} K(k,k-m, m)\;
% \overline{\widehat{f}(k)} \, \widehat{h}(k-m)\, \widehat{g}(m)\,dm\,dk
% \\
% & \quad \quad \quad=
% \int_{\R} \int_{\R} K(k,k-m, m)\;
% {\widehat{f}(-k)} \, \widehat{h}(k-m)\, \overline{\widehat{g}(-m)}\,dm\,dk
% \\
% & \quad \quad \quad=
% \int_{\R} \int_{\R} K(-m,k-m, -k)\;
% {\widehat{f}(m)} \, \widehat{h}(k-m)\, \overline{\widehat{g}(k)}\,dk\,dm
%  \\
% & \quad \quad \quad=
% \int_{\R} \int_{\R} K(-m,k-m, -k)\;
% {\widehat{f}(m)} \, \widehat{h}(k-m)\, \overline{\widehat{g}(k)}\,dm\,dk.
%\end{align*}
\qed

\medskip

\begin{corollary} \label{cor32}
Let $\varepsilon < \varepsilon_0$ and $\varepsilon_0$ be sufficiently small.
For $\ell \ge 1$, the energy $\mathcal{E}_{\ell}$ is equivalent to $ \big(\|R_{-1}\|_{H^{\ell}} +  \|R_1\|_{H^{\ell}}\big)^2$,
%  and sufficiently small $\veps>0$.
i.e. there are constants $C_1, C_2 > 0$ such that
\begin{align*}
\big( \|R_{-1}\|_{H^{\ell}} + \|R_1\|_{H^{\ell}} \big)^2
\quad \le  \quad
C_1 \,\mathcal{E}_{\ell}
 \quad \le \quad
C_2 \,
\big( \|R_{-1}\|_{H^{\ell}} +\|R_1\|_{H^{\ell}}\big)^2.
\end{align*}
\end{corollary}

\textbf{Proof.}
That  $\hat{\vartheta}(k)=1$ is true for $k$ outside of the compact set $[-\delta, \delta]$ gives the a priori estimate
\begin{align}
\label{int und theta}
\int_{\R}\partial_x^{\ell} f\, \partial_x^{m} \vartheta g\, dx = \int_{\R}\partial_x^{\ell} f\,  \partial_x^{m} g\,dx+ \mathcal{O}(\Vert f\Vert_{L^2} \Vert g\Vert_{L^2})\,.
\end{align}
With \eqref{theta^-1 k}, we also get
\begin{align}
\label{int und theta-1}
\int_{\R}\partial_x^{\ell} f\, \partial_x^{m+1} \vartheta^{-1} g\, dx = \int_{\R}\partial_x^{\ell} f\,  \partial_x^{m+1} g\, dx+\mathcal{O}(\Vert f\Vert_{L^2} \Vert g\Vert_{L^2})\,.
\end{align}
Due to \eqref{int und theta-1}, and, \eqref{Njj absch} and \eqref{Nj-j absch}, we have
\begin{align*}
E_{\ell} 
&= \sum_{j_1\in \{\pm1\}} \Big( \frac{1}{2}   \big\Vert \partial_x^{\ell} R_{j_1}) \big\Vert_{L^2}^2  +     \veps \sum_{j_2 \in \{\pm1\}} \int_{\R}\partial_x^{\ell} R_{j_1} \,\partial_x^{\ell} N_{j_1 j_2}(\psi_c,R_{j_2})\,dx \Big)
\\[2mm]
&   \quad
+ \varepsilon \,
\mathcal{O} 
\big( \|R_{-1}\|_{H^{1}}^2 + \|R_1\|_{H^{1}}^2 \big) \,.
\end{align*}
%such that we only have to look at the regularity of the terms
%\begin{align*}
%\int_{\R}\partial_x^{\ell} R_{j_1} \partial_x^{\ell} N_{j_1 j_2}(\psi_c,R_{j_2})\,dx.
%\end{align*}
For $(j_1, j_2)=(j, -j)$,   using Cauchy-Schwarz and \eqref{Nj-j absch} yields
\begin{align*}
\varepsilon \, \int_{\R}\partial_x^{\ell} R_{j} \, \partial_x^{\ell} N_{j -j}(\psi_c,R_{-j})\,dx
= \varepsilon \, \mathcal{O}(\Vert R_{-1}\Vert_{H^{\ell}} \Vert R_1 \Vert_{H^{\ell}}).
\end{align*}
For $(j_1, j_2)=(j, j)$, there is up to one  additional derivative falling on $\partial_x^{\ell} R_j$.
Using Leibniz's rule and \eqref{Njj absch}, we get
\begin{align*}
\varepsilon \,\int_{\R}\partial_x^{\ell} R_{j} \,\partial_x^{\ell} N_{j j}(\psi_c,R_{j})\,dx
&=\varepsilon \,\int_{\R}\partial_x^{\ell} R_{j} \, N_{j j}(\psi_c, \partial_x^{\ell}R_{j})\,dx+
\varepsilon \,
\mathcal{O}(\Vert R_j \Vert_{H^{\ell}}^2) \,.
\end{align*}
Lemma \ref{lemma verallgemeinerte partielle Integration} gives us
\begin{align*}
\varepsilon \,\int_{\R}\partial_x^{\ell} R_{j} \, N_{j j}(\psi_c, \partial_x^{\ell}R_{j})\,dx
= \frac{1}{2}\,
 \varepsilon \,\int_{\R}\partial_x^{\ell} R_{j} \,\big(N_{j j}(\psi_c, \partial_x^{\ell}R_{j}) + N_{j j}^*(\psi_c, \partial_x^{\ell}R_{j}) \big)\,dx \,
\end{align*}
where 
\begin{align*}
%\label{N^*}
\widehat{N}_{j_1 j_2}^* (h,f)(k) :=  \int_{\R} n_{j_1 j_2}(-m,k-m,-k) \widehat{h}(k-m) \widehat{f}(m)\,dm\,.
\end{align*}
Due to the skew-symmetry of $\rho$ and $\omega$, 
we have for $|k| \rightarrow \infty$:
\begin{align*}
&n_{j j}(k,k-m,m)+n_{j j}(-m,k-m,-k)
=
\frac{\rho(k)-\rho(m) }{ \omega(k) -\omega(m)+j \, \omega(k-m)}\, \chi_{c}(k-m) 
.
\end{align*}
Using Taylor to expand $\rho(k)$ in the point $m$ and
exploiting \eqref{as jj w>1},  \eqref{as jj w=1} and  \eqref{as jj w<1} yields
\begin{align*}
&n_{j j}(k,k-m,m)+n_{j j}(-m,k-m,-k)
\\[2mm]
%& \qquad
%=
%\frac{\rho(k)-\rho(m) }{ \omega(k) -\omega(m)+j \, \omega(k-m)}\, \chi_{c}(k-m) \,  && \quad \text{for } |k| \rightarrow \infty
%\\[2mm]
& \qquad  \qquad \qquad
=
\frac{\rho'(m) (k-m) +\mathcal{O} \big(\rho''(m) \big) }{ \omega(k) -\omega(m)+j \, \omega(k-m)}\, \chi_{c}(k-m) 
= \mathcal{O}\big(\, \chi_{c}(k-m) \big)  
&&  \text{for } |k| \rightarrow \infty  \,,
\end{align*}
due to  \eqref{rho BED2}, \eqref{ord nimmt ab rho} and \eqref{ord nimmt ab omega}. 
With  Cauchy-Schwarz, the Plancherel theorem and Young's inequality one now obtains 
\begin{align*}
\varepsilon \,\int_{\R}\partial_x^{\ell} R_{j} \, N_{j j}(\psi_c, \partial_x^{\ell}R_{j})\,dx
&=\varepsilon \,
\mathcal{O}(\Vert R_j \Vert_{H^{\ell}}^2)
%\\[2mm]
%&
%= \frac{1}{2}
%\int_{\R}\partial_x^{\ell} R_{j} \big( N_{j j}(\psi_c, \partial_x^{\ell}R_{j}) +N_{j j}^*(\psi_c, \partial_x^{\ell}R_{j})\big) \,dx 
%\\[2mm]
%& \qquad +\mathcal{O}(\Vert R_j \Vert_{H^{\ell}}^2),
\end{align*}
%We now obtained
%\begin{align*}
%{E}_{\ell}
% =
%\frac{1}{2} \big(
%\| \partial_x^{\ell} R_{-1}\|_{L^2}^2 +  \| \partial_x^{\ell} R_1\|_{L^2}^2\big)
% +
%\varepsilon \, \mathcal{O}
%\big(
%\big( \|R_{-1}\|_{H^{\ell}}+ \|R_1\|_{H^{\ell}} \big)^2 
%\big).
%\end{align*} 
such that the statement follows with lemma \ref{lemma E_0}.
\qed
\medskip

\subsection{Closing the energy estimates}
\label{SEC Closing the energy estimates}
For the Gronwall argument we are aiming for, the evolution of the energy has to be estimated against the energy itself.
More precisely, we need to estimate the evolution of $\mathcal{E}_{\ell}$ against terms involving no higher Sobolev norm  of the error  than the $H^{\ell}$-norm.
\begin{lemma}
\label{lemma d/dt E_l}
For $\ell \ge 1$, we have
\begin{align}
\label{Vell}
\partial_t E_{\ell} =
 \sum_{i=0}^4 I_i + \veps^{2}\, \mathcal{O}(\mathcal{E}_{\ell}+1),
\end{align}
where 
\begin{align*}
I_0= \varepsilon^2 \, \sum_{j_1, j_3 \in \{\pm1\}} j_1 \int_{\R} \partial^{\ell}_x R_{j_1}\, i \rho \partial_x^{\ell} \vartheta^{-1}(R_Q \vartheta R_{j_3})\,dx \,,
\end{align*}
\begin{align*}
I_1+I_2&=\varepsilon^2 \,
\sum_{j_1,j_3 \in \{\pm1\}}\!  j_1 \Big(\, \int_{\R}  i \rho \partial^{\ell}_x   \vartheta^{-1}(R_{\psi} \vartheta R_{j_3})\, \partial^{\ell}_x \vartheta^{-1}N_{j_1j_1}(\psi_c, R_{j_1})\,dx
\\[1mm] \nonumber
&
\quad
\qquad\qquad\qquad\qquad + \int_{\R} \partial^{\ell}_x R_{j_1}\, \partial^{\ell}_x \vartheta^{-1}N_{j_1j_1}\big(\psi_c, i\rho \vartheta^{-1} (R_{\psi} \vartheta R_{j_3})\big)\,dx \,\Big) \, , 
\end{align*}
\begin{align*}
I_3+I_4&= \varepsilon^2 \,
\sum_{j_1,j_3 \in \{\pm1\}}\!  j_1\Big(\, \int_{\R}  i \rho \partial^{\ell}_x   \vartheta^{-1}(R_{\psi} \vartheta R_{j_3})\, \partial^{\ell}_x \vartheta^{-1}N_{j_1-j_1}(\psi_c, R_{-j_1})\,dx
\\[1mm] \nonumber
&
\quad
\qquad\qquad\qquad\qquad - \int_{\R} \partial^{\ell}_x R_{j_1}\, \partial^{\ell}_x \vartheta^{-1}N_{j_1-j_1}\big(\psi_c, i\rho \vartheta^{-1} (R_{\psi} \vartheta R_{j_3})\big)\,dx \,\Big) \, ,
\end{align*}
\begin{align}
\label{R_Q}
R_Q= \psi_Q+ \frac{1}{2} \varepsilon^{\beta-2} (\vartheta R_{-1}+\vartheta R_{1})\,.
\end{align}
\end{lemma}

%\begin{remark}
%Due to \eqref{int und theta-1} and \eqref{theta-1rho} the term $ \varepsilon^{2} \,V_{\ell}$ indeed has the desired $\varepsilon^2$-order.
%\end{remark}

\textbf{Proof.}
%We have
%\begin{align*}
%\partial_t E_{\ell} &= \sum_{j_1 \in \{\pm1\}} \Big( \int_{\R} \partial^{\ell}_x R_{j_1}\, \partial_t\partial^{\ell}_x R_{j_1} \,dx\,
%\\[2mm] 
%&
%\quad
%+   \veps \sum_{j_1,j_2 \in \{\pm1\}} \Big( \int_{\R} \partial_t\partial^{\ell}_x R_{j_1}\, \partial^{\ell}_x \vartheta^{-1}N_{j_1j_2}(\psi_c,R_{j_2})\,dx
%\\[2mm]
% &
% \qquad
%\qquad\qquad\qquad
%+ \int_{\R} \partial^{\ell}_x R_{j_1}\, \partial^{\ell}_x \vartheta^{-1} \partial_t N_{j_1j_2}(\psi_c, R_{j_2})\,dx\,
%\Big)
%.
%\end{align*}
Using the error equations \eqref{R_j} and
exploiting 
$
R_{\psi}= \psi_c+ \varepsilon R_Q,
%\qquad \text{where} \qquad
%R_Q= \psi_Q+ \frac{1}{2} \varepsilon^{\beta-2} (\vartheta R_{-1}+\vartheta R_{1})
$
we get
\begin{align*}
\partial_t E_{\ell} 
=& 
\sum_{j_1 \in \{\pm1\}} \, j_1 \int_{\R} \partial^{\ell}_x R_{j_1}\, i\omega \partial^{\ell}_x R_{j_1} \,dx
\\[2mm] 
&
%\qquad\qquad\qquad\qquad\;\; 
+\, \veps \sum_{j_1, j_2 \in \{\pm1\}}\! \Big( \,
j_1 \int_{\R} \partial^{\ell}_x R_{j_1}\, i \rho \partial^{\ell}_x \vartheta^{-1}({\psi_c}  \vartheta R_{j_2})\,dx \; 
\\[2mm]
& \qquad\qquad\qquad \quad + j_1  \int_{\R} i \omega \partial^{\ell}_x R_{j_1}\, \partial^{\ell}_x\vartheta^{-1}N_{j_1j_2}(\psi_c, R_{j_2})\,dx
\\[2mm]
&\qquad\qquad\qquad \quad + j_2  \int_{\R} \partial^{\ell}_x R_{j_1}\, \partial^{\ell}_x \vartheta^{-1}N_{j_1j_2}(\psi_c, i \omega R_{j_2})\,dx
\\[2mm]
&\qquad\qquad\qquad \quad - \int_{\R} \partial^{\ell}_x R_{j_1}\, \partial^{\ell}_x \vartheta^{-1}N_{j_1j_2}(i \omega \psi_c, R_{j_2})\,dx
\\[2mm]
&\qquad\qquad\qquad \quad
+ \int_{\R} \partial^{\ell}_x R_{j_1}\, \partial^{\ell}_x \vartheta^{-1}N_{j_1j_2}(\partial_t \psi_c +i \omega \psi_c, R_{j_2})\,dx 
 \Big)
 \\[3mm]
&
+\,\veps^{2} \sum_{j_1, j_2 \in \{\pm1\}} j_1 \int_{\R} \partial^{\ell}_x R_{j_1}\, i \rho \partial_x^{\ell} \vartheta^{-1}(R_Q \vartheta R_{j_2})\,dx 
 \\[2mm]
&
+\, \veps^{2}  \sum_{j_1,j_2,j_3 \in \{\pm1\}}\!  \Big(\, j_1\int_{\R}  i \rho \partial^{\ell}_x   \vartheta^{-1}(R_{\psi} \vartheta R_{j_3})\, \partial^{\ell}_x \vartheta^{-1}N_{j_1j_2}(\psi_c, R_{j_2})\,dx
\\[1mm]
&\qquad\qquad\qquad\qquad +j_2 \int_{\R} \partial^{\ell}_x R_{j_1}\, \partial^{\ell}_x \vartheta^{-1}N_{j_1j_2}\big(\psi_c, i\rho \vartheta^{-1} (R_{\psi} \vartheta R_{j_3})\big)\,dx \,\Big) 
\\[3mm]
&
+ \sum_{j_1 \in \{\pm1\}}  \int_{\R} \partial^{\ell}_x R_{j_1}\, \veps^{-\beta} \partial^{\ell}_x \vartheta^{-1} \mathrm{Res}_{u_{j_1}}\!(\veps\Psi)\,dx\,
\\[2mm]
& 
+\veps \sum_{j_1, j_2 \in \{\pm1\}}\!
\Big(
 \int_{\R} \veps^{-\beta} \partial^{\ell}_x \vartheta^{-1} \mathrm{Res}_{u_{j_1}}\!(\veps\Psi)\, \partial^{\ell}_x \vartheta^{-1}N_{j_1j_2}(\psi_c,R_{j_2})\,dx
\\[2mm]& \qquad\qquad \qquad\quad 
+ \int_{\R} \partial^{\ell}_x R_{j_1}\, \partial^{\ell}_x \vartheta^{-1}N_{j_1j_2} \big(\psi_c,\veps^{-\beta} \vartheta^{-1} \mathrm{Res}_{u_{j_2}}\!(\veps\Psi)\big)\,dx\,
\Big)
\end{align*}
where 
\begin{align*}
R_Q= \psi_Q+ \frac{1}{2} \varepsilon^{\beta-2} (\vartheta R_{-1}+\vartheta R_{1}).
\end{align*}
Exploiting the skew symmetry of $i\omega$ in the third integral and then using \eqref{nf}, we get
\begin{align*}
\partial_t E_{\ell} 
=& \,
\varepsilon \! \sum_{j_1, j_2 \in \{\pm1\}}\! \Big( \,
j_1 \int_{\R} \partial^{\ell}_x R_{j_1}\, i \rho \partial^{\ell}_x \vartheta^{-1} \big({\psi_c}  (\vartheta- \vartheta_{\varepsilon, \infty}) R_{j_2} \big)\,dx \; 
\\[2mm]
& \qquad\qquad\qquad
+ \int_{\R} \partial^{\ell}_x R_{j_1}\, \partial^{\ell}_x \vartheta^{-1}N_{j_1j_2}(\partial_t \psi_c +i \omega \psi_c, R_{j_2})\,dx 
 \Big)
 \\[3mm]
&
+\sum_{i=0}^4 I_i  \,
+ \sum_{j_1 \in \{\pm1\}}  \int_{\R} \partial^{\ell}_x R_{j_1}\, \veps^{-\beta} \partial^{\ell}_x \vartheta^{-1} \mathrm{Res}_{u_{j_1}}\!(\veps\Psi)\,dx\,
\\[2mm]
& 
+\varepsilon \! \sum_{j_1, j_2 \in \{\pm1\}}\!
\Big(
 \int_{\R} \veps^{-\beta} \partial^{\ell}_x \vartheta^{-1} \mathrm{Res}_{u_{j_1}}\!(\veps\Psi)\, \partial^{\ell}_x \vartheta^{-1}N_{j_1j_2}(\psi_c,R_{j_2})\,dx
\\[2mm]& \qquad\qquad\qquad \quad
+ \int_{\R} \partial^{\ell}_x R_{j_1}\, \partial^{\ell}_x \vartheta^{-1}N_{j_1j_2} \big(\psi_c,\veps^{-\beta} \vartheta^{-1} \mathrm{Res}_{u_{j_2}}\!(\veps\Psi)\big)\,dx\,
\Big) 
\,.
\end{align*}
We now show that all rest terms can be estimated against $\veps^{2}\, \mathcal{O}(\mathcal{E}_{\ell}+1)$.
Thereby we will especially take advantage of corollary \ref{cor32} and \eqref{RES3}.
\\
Using \eqref{int und theta-1}, Cauchy-Schwarz and the fact that 
$\hat{\vartheta}- \hat{\vartheta}_{\varepsilon, \infty}=
\mathcal{O}(\varepsilon)
$ has compact support, the first integral
is estimated against $\varepsilon^{2}\, \mathcal{O}(\mathcal{E}_{\ell})$.

The second integral in the above evolution equality is $\veps^{3}\, \mathcal{O}(\mathcal{E}_{\ell})$ due to the estimate \eqref{dtpsi}. 
In order to see this we
first use \eqref{int und theta-1}, then we proceed as in the proof of \eqref{cor32} in order to estimate without  losing  regularity.

The last three integrals are $\veps^{2}\, \mathcal{O}(\mathcal{E}_{\ell}+1)$  due to \eqref{RES1}.
To see this, we use first \eqref{int und theta-1}, then integration by parts to shift some derivatives away from $ R_{\pm1}$,
and finally Cauchy-Schwarz together with \eqref{Njj absch} and \eqref{Nj-j absch}. We also exploit the estimate $\sqrt{x} \le |x| +1$ .
\qed

\medskip

The terms $I_0-I_5$  contain integrals, in which  too many additional derivatives are falling on the error functions such that
a $H^{\ell}$-energy estimate is impossible for $\deg^{*}(\rho)>0$.

In the following we will assume  $\varepsilon_0$ to be  chosen such small that 
\begin{align}
\label{epsilon E_l le 1}
\varepsilon \,\mathcal{E}_{\ell} \le 1 \,,
\end{align}
for $0<\varepsilon < \varepsilon_0 $.
%This assumption implies,  for instance 
%$\varepsilon \,\mathcal{E}_{\ell}^{3/2}= \mathcal{O}(\mathcal{E}_{\ell})$,
% $\Vert {R}_{\psi} \Vert_{H^{\ell}}= \mathcal{O}(\varepsilon^{-1/2})$ 
%and  $ \Vert {R}_{\psi} \Vert_{C^{\ell-1}}= \mathcal{O}(1)$.
This assumption is possible since there is some $T(\varepsilon) >0$ such that the $H^{\ell}$-norms of $R_{-1}(t)$ and $R_{1}(t)$ can be  uniformly bounded for $0\le t \le T(\varepsilon) $. 
When the energy estimates do close,    Gronwall's inequality  provides  $T(\varepsilon) \ge T_0 \,\varepsilon^{-2}$.

%We are now using the notation
% $ \lfloor x \rfloor = \max \{z \in \mathbb{Z}: z \le x \}$ and 
%$ \lceil \alpha \rceil := \min \{z \in \mathbb{Z}:  z \ge \alpha \}$.
\begin{lemma}
\label{artkel regloesung}
Let $0< \varepsilon< \varepsilon_0$ and $\ell \ge \lceil \deg(\omega) \rceil+ \lceil \deg^*(\rho) \rceil +1$.
\\
Let $\gamma$ be a pseudo-differential operator given by its symbol in Fourier space such that 
\begin{align}
\label{ord nimmt ab}
\deg^*(\gamma^{(l)}) \le  \deg^* (\gamma^{(l-1)})-1   \,,
\end{align} 
as long as  $\gamma^{(l)} \neq 0$.
We assume 
 $\gamma$
to be not depending on $\varepsilon$  and  to  be  either given by $\gamma=i \sigma$, where
 $\sigma \in C^{\lceil \deg^*(\sigma) \rceil}(\mathbb{R},\mathbb{R})$ is an odd function with
$ \deg^*(\sigma) \le \deg(\omega)$, or by $\gamma= \upsilon$
where $\upsilon \in C^{\lceil \deg^*(\upsilon) \rceil}(\mathbb{R},\mathbb{R})$ 
is an even function with $\deg^*(\gamma)\le \deg(\omega') $.
\\
If $\gamma=i \sigma$, 
let $f=h$ be a function with
\begin{align}
\label{h prop}
\Vert h\Vert_{H^{\lceil \deg^*(\sigma) \rceil +\lceil \deg(\omega) \rceil}}+ \Vert \partial_t  h \Vert_{H^{\lceil \deg^*(\sigma)\rceil -1}}&= \mathcal{O}(\varepsilon^{-1/2}),
% &&
\\[2mm]
\nonumber
\Vert h\Vert_{C^{\lceil \deg^*(\sigma) \rceil-1}}+\Vert \partial_t h\Vert_{C^{\lceil \deg^*(\sigma) \rceil-1}}&=\mathcal{O}(1) \,.
\end{align}
If $\gamma= \upsilon$, let $f=g$ be a function with $\Vert \partial_x^{-1} g\Vert_{\infty} = \mathcal{O}(\varepsilon^{-1})$,
\begin{align}
\label{g prop}
\Vert g\Vert_{H^{\lceil \deg^*(\upsilon) \rceil + \lceil \deg(\omega) \rceil}}+ \Vert \partial_t  g \Vert_{H^{\lceil \deg^*(\upsilon) \rceil}}&= \mathcal{O}(\varepsilon^{-1/2}),
\\[2mm]  
\nonumber
\Vert  g\Vert_{C^{\lceil \deg^*(\upsilon) \rceil-1}} +\Vert \partial_t \partial_x^{-1} g\Vert_{C^{\lceil \deg^*(\upsilon) \rceil}}&=\mathcal{O}(1).
\end{align}
Suppose  
\begin{align}
\label{f prop}
\Vert f\Vert_{H^{\lceil \deg^*(\gamma) \rceil +\lceil \deg(\omega) \rceil}}= \mathcal{O}(1)
 && \textit{or} 
 &&
 \Vert \widehat{f} \Vert_{L^1(\lceil \deg^*(\gamma) \rceil +\lceil \deg(\omega) \rceil)}= \mathcal{O}(1),
\end{align}
then there exists an $\varepsilon_0>0$ such that for all $j_1,j_2 \in\{\pm 1\}$
there is an expression $\mathcal{D}$  with 
\begin{align*}
\varepsilon^{2} \,\mathcal{D}= \varepsilon \,\mathcal{O}(\mathcal{E}_{\ell})
\end{align*}
and
\begin{align}
%\label{umform}
\varepsilon^{2}  \int_{\mathbb{R}}&  \gamma \partial_x^{\ell} R_{j_1} \, \partial_x^{\ell} R_{j_2} \, f \,dx
=\varepsilon^2 \partial_t \mathcal{D} 
+\varepsilon^{2}\,\mathcal{O}(\mathcal{E}_{\ell}+1) \,.
\end{align}
\end{lemma}
\medskip

We will will show the proof of this lemma  after  the following corollary.

\begin{corollary}
\label{corollary checkEnergy}
Let $\ell \ge \lceil \deg(\omega) \rceil + \lceil \deg^*(\rho) \rceil +1$.
\\
For  $\varepsilon_0$ sufficiently  small and $0<\varepsilon < \varepsilon_0$, there exists an energy $\tilde{\mathcal{E}}_{\ell}$  and some constants $c, C >0$ such that 
\begin{align}
\big( \Vert R_{-1} \Vert_{H^{\ell}} +\Vert R_{1} \Vert_{H^{\ell}} \big)^2 \le c\, \tilde{\mathcal{E}}_{\ell} \le C \, \big( \Vert R_{-1} \Vert_{H^{\ell}} +\Vert R_{1} \Vert_{H^{\ell}} \big)^2
\end{align}
and
\begin{align*}
\partial_t \tilde{\mathcal{E}}_{\ell} \le \varepsilon^2 \,\mathcal{O} \big( \tilde{\mathcal{E}}_{\ell}+1 \big) 
\,.
\end{align*}

\end{corollary}

\textbf{Proof.}
According to  lemma \ref{lemma d/dt E_0} and lemma \ref{lemma d/dt E_l}, 
we have
\begin{eqnarray*}
\partial_t {\mathcal{E}}_{\ell}  &=&
%\varepsilon^2 \, V_{\ell}
%+\, \veps^{2}\, \mathcal{O}(\mathcal{E}_{\ell}+1),
%\\[2mm] 
% &=& \veps^{2} \sum_{j_1,j_3  \in \{\pm1\}} \Big(\,
%j_1 \int_{\R} \partial^{\ell}_x R_{j_1}\, i \rho \partial^{\ell}_x \vartheta^{-1}({ R_Q} \vartheta R_{j_3})\,dx \; 
%\\[2mm] 
%&&\qquad\qquad\qquad+j_1 \int_{\R}  i \rho \partial^{\ell}_x \vartheta^{-1}(R_{\psi}\vartheta R_{j_3})\, \partial_x^{\ell}\vartheta^{-1}N_{j_1j_1}(\psi_c, R_{j_1})\,dx\\[2mm]
%&&\qquad\qquad\qquad+ j_1 \int_{\R} \partial^{\ell}_x R_{j_1}\, \partial^{\ell}_x  \vartheta^{-1}N_{j_1j_1}(\psi_c, i\rho \vartheta^{-1}(R_{\psi} \vartheta R_{j_3}))\,dx 
%\\[2mm]&&
%\qquad\qquad\qquad+j_1 \int_{\R}  i \rho \partial^{\ell}_x \vartheta^{-1}(R_{\psi} \vartheta R_{j_3})\, \partial_x^{\ell}\vartheta^{-1}N_{j_1-j_1}(\psi_c, R_{-j_1})\,dx\; \\[2mm]
%&&\qquad\qquad\qquad- j_1 \int_{\R} \partial^{\ell}_x R_{j_1}\, \partial^{\ell}_x\vartheta^{-1}N_{j_1-j_1}(\psi_c, i\rho \vartheta^{-1}(R_{\psi} \vartheta R_{j_3}))\,dx\, \Big) \,
% \\[2mm]
%&&\quad +\, \veps^{2}\, \mathcal{O}(\mathcal{E}_{\ell}+1)\\[2mm] 
%&=:&
 \sum_{i=0}^4 I_i + \veps^{2}\,\mathcal{O}(\mathcal{E}_{\ell}+1)\,.
\end{eqnarray*}
First, we analyze the term $I_0$. 
\\
To easier keep track of the terms containing the highest derivatives of the error, we  introduce the notations
\begin{align}
\label{def tildeRpsi}
\tilde{R}_{\psi}:= \psi+ \varepsilon^{\beta-1} \vartheta (R_{1}+R_{-1}),  \qquad 
\tilde{R}_Q:= \psi_Q+ \varepsilon^{\beta-2} \vartheta (R_{-1}+ R_{1})\,.
\end{align}
For $N \in \mathbb{N}$ and $ \ell \ge 2N +1$, Leibniz's rule yields
\begin{align}
\label{LeibnizR}
\partial_x^{\ell} \big( R_{\psi} 
\vartheta (R_{1}+R_{-1}) \big)
 & 
 =
\sum_{n=0}^{N} \binom{\ell}{n}
 \partial_x^{n} \tilde{R}_{\psi} \, \partial_x^{\ell-n}  \vartheta (R_{1}+R_{-1})  
 \\ \nonumber
 & \quad
+
\sum_{n=N+1}^{ \ell-N-1} \binom{\ell}{n}
 \partial_x^{n} R_{\psi} \, \partial_x^{\ell-n}  \vartheta (R_{1}+R_{-1})
  \\ \nonumber
 & \quad 
  +
\sum_{n=\ell- N}^{\ell} \binom{\ell}{n}
    \partial_x^{n} \psi\, \partial_x^{\ell-n}  \vartheta (R_{1}+R_{-1})  \,.
\end{align}
So only the first term is important when Leibniz's rule is applied.

Using \eqref{int und theta-1}, the skew symmetry of $i \rho$, Leibniz's rule and \eqref{int und theta}, we get
\begin{align}
\label{I_0}
I_0 &=\veps^{2} \sum_{j_1,j_3  \in \{\pm1\}} 
j_1 \int_{\R} \partial^{\ell}_x R_{j_1}\, i \rho \partial^{\ell}_x \vartheta^{-1}({ R_Q} \vartheta R_{j_3})\,dx
\\ \nonumber
&=
\veps^{2} \sum_{j_1,j_3  \in \{\pm1\}} 
  \sum_{n=0}^{N} \binom{\ell}{n}
\int_{\R}  i \rho  \partial^{\ell}_x R_{j_1}\, 
  \partial^{ \ell-n}_x R_{j_3} \, \partial^{n}_x \tilde{R}_Q  \,dx
+\, \veps^{2}\,\mathcal{O}(\mathcal{E}_{\ell}+1)\,
\, ,
\end{align}
where $N:= \lceil \deg^*(\rho \rceil -1$.
With integration by parts and lemma \ref{artkel regloesung}, one obtains
\begin{align*}
I_0 
&= \varepsilon^2 \, \partial_t  {\mathcal{D}}_{0}+ \veps^{2}\,\mathcal{O}(\mathcal{E}_{\ell}+1)\,
\end{align*}
for some ${\mathcal{D}}_{0}$ with $\varepsilon^2 {\mathcal{D}}_{0}=  \varepsilon\,\mathcal{O}(\mathcal{E}_{\ell})$.

\medskip

Now, we analyze the term $I_1+I_2$. 
\\
Using \eqref{int und theta-1} and \eqref{Njj absch}, we get
\begin{align}
\label{I_1+I_2}
I_1+ I_2 
&= 
\veps^{2} \sum_{j_1,j_3  \in \{\pm1\}} 
j_1
\Big(\,
 \int_{\R}  i \rho \partial^{\ell}_x \vartheta^{-1}(R_{\psi}\vartheta R_{j_3})\, \partial_x^{\ell}N_{j_1j_1}(\psi_c, R_{j_1})\,dx
\\ \nonumber
& \quad\quad\quad \quad \quad \quad \quad
+ \int_{\R} \partial^{\ell}_x R_{j_1}\, \partial^{\ell}_x  N_{j_1j_1}
\big(\psi_c, i\rho \vartheta^{-1}(R_{\psi} \vartheta R_{j_3}) \big)\,dx 
\Big)
 \\ \nonumber
& \qquad   +\, \veps^{2}\,\mathcal{O}(\mathcal{E}_{\ell}+ \veps^{\beta -1} \mathcal{E}^{3/2}_{\ell})\,. 
\end{align}
Due to \eqref{Njj absch}, applying Leibniz's rule gives
\begin{align*}
I_1+ I_2 
&= \veps^{2} \sum_{j_1, j_3 \in \{\pm1\}} 
j_1
\Big(\,\int_{\R}  i \rho \partial^{\ell}_x \vartheta^{-1}(R_{\psi} \vartheta R_{j_3})\, N_{j_1j_1}(\psi_c, \partial^{\ell}_x R_{j_1})\,dx
%\\
%&\qquad\qquad\qquad \quad
%+ \ell   \int_{\R}  i \rho \partial^{\ell}_x \vartheta^{-1}(R_{\psi} \vartheta R_{j_3})\, N_{j_1j_1}(\partial_x\psi_c, \partial^{\ell-1}_x R_{j_1})\,dx
\\
&\qquad\qquad\qquad 
+\!  \sum_{m=1}^{\lceil \deg^*(\rho)\rceil}\! \binom{\ell}{m} \!
   \int_{\R}  i \rho \partial^{\ell}_x \vartheta^{-1}(R_{\psi} \vartheta R_{j_3})\, N_{j_1j_1}(\partial_x^m \psi_c, \partial^{\ell-m}_x R_{j_1})\,dx
\\
&\qquad\qquad\qquad + \int_{\R} \partial^{\ell}_x R_{j_1} \,N_{j_1j_1}\big(\psi_c, i \rho \partial^{\ell}_x \vartheta^{-1}(R_{\psi} \vartheta R_{j_3})\big)\,dx
%\\
%&\qquad\qquad\qquad \quad+  \ell   \int_{\R} \partial^{\ell}_x R_{j_1} \,N_{j_1j_1}\big(\partial_x \psi_c, i \rho \partial^{\ell-1}_x \vartheta^{-1}(R_{\psi} \vartheta R_{j_3})\big)\,dx
\\
&\qquad\qquad\qquad 
+  \!\sum_{m=1}^{\lceil \deg^*(\rho)\rceil}\! \binom{\ell}{m} \!
\int_{\R} \partial^{\ell}_x R_{j_1} \,N_{j_1j_1}\big(\partial_x^m \psi_c, i \rho \partial^{\ell-m}_x \vartheta^{-1}(R_{\psi} \vartheta R_{j_3})\big)\,dx
\Big)
 \\[2mm]
&\qquad \qquad+\, \veps^{2}\,\mathcal{O}(\mathcal{E}_{\ell}+ \veps^{\beta -1} \mathcal{E}^{3/2}_{\ell})\,.
\end{align*}
By  using lemma \ref{lemma verallgemeinerte partielle Integration}, the skew symmetry of $i \rho$  and  \eqref{int und theta-1} we get
\begin{align*}
I_1+ I_2 &
= 
\varepsilon^{2}  \sum_{j_1, j_3 \in \{\pm1\}} j_1 \Big(\,  \int_{\R}  \partial^{\ell}_x (R_{\psi} \vartheta R_{j_3})\, 
i \rho \big(
N_{j_1j_1}(\psi_c, \partial^{\ell}_x R_{j_1})
+ \,N_{j_1j_1}^*(\psi_c,\partial^{\ell}_x R_{j_1} )
\big)
\,dx
\\
&\qquad\qquad\qquad \qquad
+ \! \sum_{m=1}^{\lceil \deg^*(\rho)\rceil}\! \binom{\ell}{m}  \!
   \int_{\R}   \partial^{\ell}_x (R_{\psi} \vartheta R_{j_3})\, i \rho  N_{j_1j_1}(\partial_x^m \psi_c, \partial^{\ell-m}_x R_{j_1})\,dx
\\
&\qquad\qquad\qquad  \qquad
+\!  \sum_{m=1}^{\lceil \deg^*(\rho)\rceil}\! \binom{\ell}{m}  \!
\int_{\R}   \partial^{\ell-m}_x (R_{\psi} \vartheta R_{j_3}) \, i \rho N_{j_1j_1}^*(\partial_x^m \psi_c,\partial^{\ell}_x R_{j_1}  )\,dx
\,
\Big)
 \\[2mm]
&\qquad +\, \veps^{2}\,\mathcal{O}(\mathcal{E}_{\ell}+ \veps^{\beta -1} \mathcal{E}^{3/2}_{\ell})\,,
\end{align*}
where
\begin{align*}
%\label{N^*}
\widehat{N}_{j_1 j_2}^* (\psi_c,f)(k) :=  \int_{\R} n_{j_1 j_2}(-m,k-m,-k) \widehat{\psi}_c(k-m) \widehat{f}(m)\,dm\,.
\end{align*}
If we now look at
\begin{align*}
 i \rho(k) \,\big(
n_{j j}&(k,k-m,m)+n_{j j}(-m,k-m,-k) \big)
\\[2mm]
&= i \rho(k) \, (\rho(k)- \rho(m)) \,\frac{\chi_c(k-m)}{ \omega(k) -\omega(m)+j \, \omega(k-m)}
&& ( \text{for }|k| \rightarrow \infty) 
\,
\end{align*}
and use Taylor's theorem, the same cancellation as in the proof of corollary \ref{cor32} occurs.
By now exploiting \eqref{as jj w>1} or respectively \eqref{as jj w=1} or \eqref{as jj w<1} (with integration by parts for the third term), we get
\begin{align*}
I_1+ I_2 &
= 
\veps^{2} \sum_{j_1, j_3 \in \{\pm1\}} j_1 \, \sum_{n=1}^N
% \Big(\,
 \int_{\R}   \partial^{\ell}_x (R_{\psi}  \vartheta  R_{j_3})\, 
\beta _n \psi_c  \, \alpha_n \partial^{\ell}_x R_{j_1}
 \,dx
+\veps^{2}\,\mathcal{O}(\mathcal{E}_{\ell}+ 1)
\,,
\end{align*}
for some $N \in \mathbb{N}$ and some pseudo-differential operators $\beta_n$ and $\alpha_n$, where $\alpha_n$ is either skew-symmetric with $\deg^*(\alpha_n) \le \deg^*(\rho)$
or symmetric with  $\deg^*(\alpha_n) \le \deg^*(\rho)-1$.
With the help
of \eqref{LeibnizR}, \eqref{int und theta} and integration by parts
we can now apply lemma \ref{artkel regloesung} to obtain
\begin{align*}
I_1+ I_2 &
=  \varepsilon^2 \, \partial_t  {\mathcal{D}}_{1,2}+\veps^{2}\,\mathcal{O}(\mathcal{E}_{\ell}+ 1)\,
\end{align*}
for some ${\mathcal{D}}_{1,2}$ with $\varepsilon^2 \, {\mathcal{D}}_{1,2}= \varepsilon \,\mathcal{O}(\mathcal{E}_{\ell})$.
To apply  lemma \ref{artkel regloesung} one splits 
\begin{align*}
\partial^{m_1}_x \tilde{R}_{\psi}  \partial^{m_2}_x\beta _n \psi_c &=
\partial^{m_1}_x (\psi_c+ \varepsilon \tilde{R}_{Q})  \partial^{m_2}_x\beta _n \psi_c
\\
&
=\partial^{m_1}_x \psi_c  \partial^{m_2}_x\beta _n \psi_c
+\partial^{m_1}_x \varepsilon \tilde{R}_{Q}  \partial^{m_2}_x\beta _n \psi_c =:f_1+f_2
\end{align*}
such that $\Vert f_1 \Vert_{L^1(p)}= \mathcal{O}(1)$ and $\Vert f_2 \Vert_{H^p}= \mathcal{O}(1)$.
The estimate $\Vert \partial_x^{-1} f_1 \Vert_{\infty}= \mathcal{O}(\varepsilon^{-1})$
is obtained by 
\begin{align*}
\Vert \partial_x^{-1} f_1 \Vert_{\infty}
&=
\Vert \partial_x^{-1} ( \partial^{m_1}_x \psi_c  \partial^{m_2}_x\beta _n \psi_c )\Vert_{\infty}
\le
\int_{\mathbb{R}} |\partial^{m_1}_x \psi_c  \partial^{m_2}_x\beta _n \psi_c| \,dx
\\
&
 \le 
\Vert \partial^{m_1}_x \psi_c  \Vert_{L^2} \Vert \partial^{m_2}_x\beta _n \psi_c\Vert_{L^2}
\,.
\end{align*}
To obtain the estimate $\Vert \partial_t \partial_x^{-1} f_1 \Vert_{\infty}= \mathcal{O}(1)$ one has to proceed more carefully. To obtain this estimate one
splits $\psi_c = \psi_{-1}+\psi_{1}$ with $\psi_{\pm1}$ as in \eqref{ansatz2}. Then one exploits the fact that  the products $\psi_j \psi_j$ are 
strictly concentrated around $k= \pm 2 k_0$ in Fourier space 
and  $\partial_t (\psi_j \psi_{-j})=  \mathcal{O}(\varepsilon)$  
such that
$\Vert  \mathcal{F} \big[ \partial_t \partial_x^{-1} (\psi_j \psi_j) \big] \Vert_{L^1(p)}= \mathcal{O}(1)$ and
$\Vert \partial_t \partial_x^{-1} (\psi_j \psi_{-j}) \Vert_{\infty}= \mathcal{O}(1)$ can be obtained.
The other  estimates  are straightforward.

\medskip
Now, we analyze the term $I_3+I_4$. 
\\
Using \eqref{int und theta-1} and \eqref{Nj-j absch}, we have
\begin{align}
\label{I_3+I_4}
I_3+I_4
&=
\veps^{2} \sum_{j_1, j_3 \in \{\pm1\}}
 j_1 \Big(
 \int_{\R}  i \rho \partial^{\ell}_x \vartheta^{-1}(R_{\psi} \vartheta R_{j_3})\, \partial_x^{\ell}N_{j_1-j_1}(\psi_c, R_{-j_1})\,dx\; 
\\[2mm] \nonumber
&\qquad\qquad\qquad \quad
-  \int_{\R} \partial^{\ell}_x R_{j_1}\, \partial^{\ell}_xN_{j_1-j_1} \big(\psi_c, i\rho \vartheta^{-1}(R_{\psi} \vartheta R_{j_3})\big)\,dx\, \Big) 
\\[2mm] \nonumber
& \qquad \qquad
+ \veps^{2}\,\mathcal{O}(\mathcal{E}_{\ell}+\veps^{\beta-1} \mathcal{E}^{3/2}_{\ell}).
\end{align}
According to lemma \ref{lemma Normalformtransfo}
the $N_{j_1-j_1}(\partial_x^m \psi_c, \cdot)$  always  map $L^2 (\mathbb{R})$ on $L^2 (\mathbb{R})$.
With the help of \eqref{asyj-j w al} we can thus  proceed as before for $I_1+I_2$ after the cancellation was achieved. We apply lemma \ref{artkel regloesung}  and obtain
\begin{align*}
I_3+ I_4 
&=  \varepsilon^2 \, \partial_t \mathcal{D}_{3,4}
+ \veps^{2}\,\mathcal{O}(\mathcal{E}_{\ell}+1),
\end{align*} 
for some  $\mathcal{D}_{3,4}$  with $\varepsilon^2 \, \mathcal{D}_{3,4}=\varepsilon \, \mathcal{O}(\mathcal{E}_{\ell})$.

Choosing $\varepsilon_0$ small enough and summing up the results for $I_0$-$I_4$, we can define a modified energy 
\begin{align*}
\tilde{\mathcal{E}_{\ell}} := \mathcal{E}_{\ell} -\varepsilon^2( \mathcal{D}_{0}+ \tilde{\mathcal{D}}_{1,2}+ \mathcal{D}_{3,4})\,,
\end{align*}
%with
%\begin{align*}
%\varepsilon^2 \,(\mathcal{D}_{0}+  \tilde{\mathcal{D}}_{1,2}+  \mathcal{D}_{3,4} )= \varepsilon\, \mathcal{O} ({\mathcal{E}_{\ell}} )
%\end{align*}
such that
\begin{align*}
\partial_t \tilde{\mathcal{E}_{\ell}}\, \lesssim\, \veps^{2} \big(1+\mathcal{E}_{\ell} \big)
%\, \lesssim\, \veps^{2} \big(1+\tilde{\mathcal{E}_{\ell}} \big)
\,.
\end{align*}
Since $\tilde{\mathcal{E}_{\ell}} = \mathcal{E}_{\ell} +\varepsilon \,\mathcal{O}(\mathcal{E}_{\ell} )$, the statement   follows
with corollary \ref{cor32}.
\qed
\medskip

For the proof of lemma \ref{artkel regloesung} we use the notation
$
[\gamma,f]g:=\gamma (fg)- f \, \gamma g \, 
$  for an operator $\gamma$ and  functions $g$ and $f$.
Further, we need the following lemma.
\begin{lemma}
\label{lemma komtaylor}
Let $n \in \mathbb{N}$, and $\gamma$ be a function of $C^{n+1}(\mathbb{R})$ with $\deg^*(\gamma) \in \mathbb{R}$ for which 
\begin{align}
%\label{ord nimmt ab}
\deg^*(\gamma^{(l)}) \le  \deg^* (\gamma^{(l-1)})-1  \qquad \textit{for all }\, 1 \le l \le n+1.
\end{align}
Moreover let the operators $\gamma$ and $i^l \gamma^{(l)}$ be given by their symbols in Fourier space.
\\
Then we  have for $f ,g \in C^{\infty}_c(\mathbb{R})$:
\begin{align}
\label{komtaylor}
\big[ \gamma, g \big] f
=
\sum_{l=1}^{n} \frac{(-1)^l}{l!} \, \partial_x^l g \, i^l \gamma^{(l)}f + \mathcal{R}(f,g).
\end{align}
For  the rest-term $\mathcal{R}(f,g)$, given through
\begin{align*}
\widehat{\mathcal{R}(f,g)}
= \int_{\mathbb{R}} 
\Big(
\frac{(\cdot-m)^{n+1}}{n!} 
\int_0^1 \gamma^{(n+1)} \big(m+(\cdot-m)x \big) \,(1-x) \, dx \Big) \,
 \widehat{g}(\cdot-m) \widehat{f}(m) \, dm \,,
\end{align*}
we have the estimate
\begin{align}
\label{komrestneu}
\Vert \mathcal{R}(f,g) \Vert_{L^2}  \le \mathcal{O}(1) \,
\big\Vert  \mathcal{F}^{-1} \big[  \vert (1+|\cdot|^2)^{p/2} \,\widehat{\partial_x^{n+1}g}(\cdot)\vert \big] \big\Vert_{\infty} \Vert f\Vert_{H^p}
\end{align}
with 
$
p=\max \{\deg^*(\gamma)-n-1,\, 0\} .
$
\end{lemma}
\begin{remark}
Estimate \eqref{komrestneu}  implies 
$\Vert \mathcal{R}(f,g) \Vert_{L^2}  \le \mathcal{O}(\Vert \widehat{\partial_x^{n+1}g} \Vert_{L^1(p)}) \Vert f\Vert_{H^p} $
and, with Sobolev's embedding theorem
$\Vert \mathcal{R}(f,g) \Vert_{L^2}  \le \mathcal{O}(\Vert g\Vert_{H^{p+q}}) \Vert f\Vert_{H^p} $
for $q > 1/2$.
\end{remark}
\textbf{Proof.}
We have
\begin{align*}
\widehat{\big[ \gamma, g \big] f}= \widehat{ \gamma (g f)}-\widehat{g  \gamma  f}= \int_{\mathbb{R}} \big(\gamma(\cdot) - \gamma(m) \big) \widehat{g}(\cdot-m) \widehat{f}(m) \, dm.
\end{align*}
Using Taylor, we get
\begin{align*}
\gamma(k)- \gamma(m)&= \sum_{l=1}^{n} \frac{(k-m)^l}{l!} \gamma^{(l)}(m)+r(k,k-m,m)
\\
&=\sum_{l=1}^{n} \frac{i^l(k-m)^l}{l!} (-i)^l\gamma^{(l)}(m)+r(k,k-m,m)
,
\end{align*}
where
\begin{align*}
r(k,k-m,m)&= \frac{(k-m)^{n+1}}{n!} 
\int_0^1 \gamma^{(n+1)} \big(m+(k-m)x \big)(1-x) \, dx
\\[2mm]
&
\le  \frac{(k-m)^{n+1}}{n!}  \,\max_{x \in [0,1]} \gamma^{(n+1)} \big(m+(k-m)x \big)
\\[2mm]
&
\le
\mathcal{O}(  |k-m|^{n+1}) \big(1+ (1+|k-m|)^{\deg^*(\gamma)-n-1} \! +(1+|m|)^{\deg^*(\gamma)-n-1}  \big) .
\end{align*}
We now get
\begin{align*}
\Vert \mathcal{R}(f,g) \Vert_{L^2}
%&= \Big\Vert 
%\int_{\mathbb{R}} r(k,k-m,m)\, \widehat{g}(k-m) \,\widehat{f}(m)  \,dm
% \Big\Vert_{L^2}
%\\
&\le \Big\Vert 
\int_{\mathbb{R}} \vert r(k,k-m,m) \,  \widehat{g}(k-m)\, \widehat{f}(m) \vert \,dm
 \Big\Vert_{L^2}
 \\
&\le
\mathcal{O}(1)\,
\Big\Vert 
\int_{\mathbb{R}}  \vert (1+|k-m|^2)^{p/2} \,\widehat{\partial_x^{n+1}g}(k-m) \, (1+|m|^2)^{p/2} \,\widehat{f}(m) \vert \,dm
 \Big\Vert_{L^2} \,,
\end{align*}
with
$p=\max \{\deg^*(\gamma)-n-1,\, 0\}$.
With Plancherel's theorem, we obtain
\begin{align*}
\Vert \mathcal{R}(f,g) \Vert_{L^2} 
&\le
\mathcal{O}(1)\,
\Big\Vert \, \mathcal{F}^{-1} \big[  \vert (1+|\cdot|^2)^{p/2} \,\widehat{\partial_x^{n+1} g}(\cdot) \vert\big] \,\, \mathcal{F}^{-1} \big[\vert (1+|\cdot|^2)^{p/2} \,\widehat{f}(\cdot) \vert  \big]
 \, \Big\Vert_{L^2}
\\[2mm]
& \le
\mathcal{O}(1)\,
\big\Vert  \mathcal{F}^{-1} \big[  \vert (1+|\cdot|^2)^{p/2} \,\widehat{\partial_x^{n+1}g}(\cdot)\vert \big] \big\Vert_{\infty} \,
\Vert f \Vert_{H^p}.
\end{align*}
\qed

\medskip

\textbf{Proof of lemma \ref{artkel regloesung}.}
If $\deg^*(\gamma) \le 0$, the lemma is trivially true.
So we will in the following assume  $\deg^*(\gamma) > 0$.
Since $\deg(\omega)\ge \deg^*(\gamma) >0$, there exist some constants $D_{\omega},d_{\omega} > 0 $
such that
$
%\label{NST omega}
\vert \omega(k) \vert \ge d_{\omega}>0 $ for  $|k| \ge D_{\omega}$.
For $\gamma= \upsilon$, we can on top of that find $D_{\omega},d_{\omega} > 0 $ such that
$
%\label{NST omega'}
\vert \omega'(k) \vert \ge d_{\omega}>0 $ for  $|k| \ge D_{\omega}$
due to $\deg(\omega')\ge \deg^*(\upsilon) >0$.
\\
There is some $D \ge D_{\omega}$ and  some function $\tilde{\gamma} \in C^{\lceil \deg^*(\gamma) \rceil}(\mathbb{R},\mathbb{R})$ with \eqref{ord nimmt ab} such that 
$
\tilde{\gamma}(k)= \gamma(k)$ for $ |k|\ge D $
and
$
\tilde{\gamma}(k)= 0 $ for $ |k|\le D_{\omega}$.
Since 
\begin{align*}
\varepsilon^{2}  \int_{\mathbb{R}} \gamma \partial_x^{\ell} R_{j_1} \, \partial_x^{\ell} R_{j_2} \, f \,dx
&=
\varepsilon^{2}  \int_{\mathbb{R}} \tilde{\gamma} \partial_x^{\ell} R_{j_1} \, \partial_x^{\ell} R_{j_2} \, f \,dx
+
\varepsilon^{2}  \int_{\mathbb{R}} (\gamma-\tilde{\gamma}) \partial_x^{\ell} R_{j_1} \, \partial_x^{\ell} R_{j_2} \, f \,dx
\\[2mm]
&=
\varepsilon^{2}  \int_{\mathbb{R}} \tilde{\gamma} \partial_x^{\ell} R_{j_1} \, \partial_x^{\ell} R_{j_2} \, f \,dx
+\varepsilon^{2} \,\mathcal{O}(\mathcal{E}_{\ell}+1),
\end{align*}
we can  in the following assume that we have 
%$\gamma \in C^{\lceil \deg^*(\gamma) \rceil}(\mathbb{R},\mathbb{R})$ with \eqref{ord nimmt ab}
%and 
$\gamma(k)=0$ for $|k| \le D_{\omega}$.
This makes the operators $\frac{\gamma}{\omega}$ and  $\frac{\gamma}{\omega'}$    well-defined.
\\
As a first step, we show the following statement:
\\
\textit{
There
is an expression $\mathcal{D}$  with 
\begin{align*}
\varepsilon^2 \, \mathcal{D}= \varepsilon \, \mathcal{O}(\mathcal{E}_{\ell}),
\end{align*}
such that 
\begin{align}
\label{umform}
\varepsilon^{2}  \int_{\mathbb{R}}&  \gamma \partial_x^{\ell} R_{j_1} \, \partial_x^{\ell} R_{j_2} \, f \,dx
\\[2mm] \nonumber
&
=\varepsilon^2 \partial_t \mathcal{D} 
+\varepsilon^2 \sum_{k=1}^{\lceil\deg^*(\gamma)\rceil-1}
\int_{\mathbb{R}}  {\varsigma}_k \partial_x^{\ell} R_{j_1} \, \partial_x^{\ell} R_{j_2} \, \partial_x^k f \,dx
\\[2mm] \nonumber
& \quad
+\varepsilon^{2} \sum_{k=1}^{m_{\gamma}}
\int_{\mathbb{R}}  \gamma_k \partial_x^{\ell} R_{p_k} \, \partial_x^{\ell} R_{q_k} \, f_k \,dx
+\varepsilon^{2}\,\mathcal{O}(\mathcal{E}_{\ell}+1),
\end{align}
where the  $\varsigma_k$ and $\gamma_k$ are skew symmetric 
or symmetric operators independent of $\varepsilon$ and given by their symbol in Fourier space,
 $m_{\gamma}=m_{\gamma} \big(\deg^*(\gamma) \big) \in \mathbb{N}$,  $f_k$ are some functions and $p_k, q_k \in \{-1, 1\}$.
The functions $\varsigma_k \in C^{\lceil \deg^*(\varsigma_k )\rceil}(\mathbb{R},\mathbb{R})$ and ${\gamma}_k \in C^{\lceil \deg^*({\gamma}_k) \rceil}(\mathbb{R},\mathbb{R})$
share the property \eqref{ord nimmt ab}.
We have 
\begin{align}
\label{umform ord'}
\deg^*({\varsigma}_k) \le \deg^*(\gamma)- k.
\end{align}
If a $\gamma_k$ is a skew symmetric then
$\deg^*(\gamma_k) \le \deg(\omega)$.  If a $\gamma_k$ is symmetric then
$\deg^*(\gamma_k) \le \deg(\omega')$. 
Furthermore:
\begin{align}
\label{gutes fk}
\Big\Vert  \mathcal{F}^{-1} \big[  \vert (1+|\cdot|^2)^{p/2} \,\widehat{\partial_x^{n+1}f_k}(\cdot)\vert \big] \Big\Vert_{\infty}
= \mathcal{O}(1) \,,
\end{align}
\begin{align}
\label{umform L^2}
&\Vert f_k\Vert_{H^{\lceil \deg^*(\gamma_k) \rceil +\lceil \deg(\omega) \rceil}}+ \Vert \partial_t f_k \Vert_{H^{\lceil \deg^*(\gamma_k) \rceil }}
\\[2mm] \nonumber
&
\qquad\qquad\qquad\qquad  \le 
\varepsilon \, C_1 \,\big(\Vert f\Vert_{H^{\lceil \deg^*(\gamma_k) \rceil +\lceil \deg(\omega) \rceil}}
+ \Vert \partial_t f \Vert_{H^{\lceil \deg^*(\gamma_k) \rceil }} \big)
\\[1mm] \nonumber
& \qquad\qquad\qquad\qquad \quad
+
\varepsilon^{1/2}  \,C_2 \,
\big( \Vert \partial_x^{-1} f\Vert_{C^{\lceil \deg^*(\gamma_k) \rceil }} +\Vert \partial_t \partial_x^{-1} f\Vert_{C^{\lceil \deg^*(\gamma_k) \rceil }}
\big),
\\[3mm]
\label{umform L^infty}
&\Vert \partial_x^{-1} f_k\Vert_{\infty} +\Vert \partial_t \partial_x^{-1} f_k\Vert_{C^{\lceil \deg^*(\gamma_k) \rceil}}
\\[2mm] \nonumber
&
\qquad
\le
\varepsilon^{1/2}   \,C_1 \, \big(\Vert f\Vert_{L^2}
+ \Vert \partial_t f \Vert_{L^2} \big)
+
\varepsilon   \,C_2 \,
\Big(
\Vert \partial_x^{-1} f\Vert_{C^{\lceil \deg^*(\gamma_k) \rceil }}
+\Vert \partial_t \partial_x^{-1}f \Vert_{C^{\lceil \deg^*(\gamma_k) \rceil}} 
\big)
\, ,
\end{align}
where the constants $C_1, C_2$ depend on $\tilde{R}_{\psi}, f, \gamma$ but are independent of $\varepsilon$. 
We set $C_2:=0$, when $\gamma$ is skew symmetric i.e. $f=h$.}

\medskip
\textbf{a)}
Handling integrals of the form
\begin{align}
\label{int rjr-j beliebig}
\varepsilon^{2}  &\int_{\mathbb{R}}  \gamma \partial_x^{\ell} R_{j} \, \partial_x^{\ell} R_{-j} \, f \,dx \, .
\end{align}
%with \eqref{h prop} and \eqref{f prop}, or respectively with \eqref{g prop} and \eqref{f prop}.
%\\
By exploiting the skew symmetry of $i \omega$ and \eqref{R_j}, we have
\begin{align*}
%\label{um1}
%
&
\varepsilon^{2}  \int_{\mathbb{R}}  \gamma \partial_x^{\ell} R_{j} \, \partial_x^{\ell} R_{-j} \, f \,dx
\\[2mm]  \nonumber
& \qquad
=
\frac{1}{2}j \, \varepsilon^{2} \, \partial_t \int_{\mathbb{R}}  \frac{\gamma}{i \omega} \partial_x^{\ell} R_{j} \, \partial_x^{\ell} R_{-j} \,f \, dx
\\[2mm]  \nonumber
& \qquad
 \quad 
-\frac{1}{2}  \varepsilon^{2}  \int_{\mathbb{R}} \big[i \omega, f\big]  \frac{\gamma}{i \omega} \partial_x^{\ell} R_{j} \, \partial_x^{\ell} R_{-j} \, dx
\\[2mm]  \nonumber
& \qquad
 \quad 
-
\frac{1}{2} \varepsilon^{3} \,  \int_{\mathbb{R}}  \frac{\gamma}{i \omega} i \rho   \partial_x^{\ell} \vartheta^{-1}\big( R_{\psi} \vartheta (R_{1}+R_{-1}) \big) \, \partial_x^{\ell} R_{-j} \,f \, dx
\\[2mm]  \nonumber
& \qquad
 \quad 
+
\frac{1}{2} \varepsilon^{3} \,  \int_{\mathbb{R}}  \frac{\gamma}{i \omega} \partial_x^{\ell} R_{j} \,  i \rho   \partial_x^{\ell} \vartheta^{-1}\big( R_{\psi} \vartheta (R_{1}+R_{-1}) \big) \,f \, dx
\\[2mm]  \nonumber
& \qquad
 \quad 
-
\frac{1}{2}j \, \varepsilon^{2} \,  \int_{\mathbb{R}}  \frac{\gamma}{i \omega} \partial_x^{\ell} R_{j} \, \partial_x^{\ell} R_{-j} \, \partial_t f \, dx
\\[2mm]  \nonumber
&\qquad
\quad
-\frac{1}{2}j \, \varepsilon^{2 -\beta} \,
\int_{\mathbb{R}}  \gamma \partial_x^{\ell}  \vartheta^{-1}  {\rm Res}_{u_{j}}(\veps \Psi) \, \partial_x^{\ell} R_{-j} \, f \,dx
\\[2mm]  \nonumber
&\qquad
\quad
-\frac{1}{2}j \, \varepsilon^{2 -\beta} \,
\int_{\mathbb{R}}  \gamma \partial_x^{\ell} R_{j} \, \partial_x^{\ell}  \vartheta^{-1}  {\rm Res}_{u_{-j}}(\veps \Psi) \, f \,dx
. 
\end{align*}
The first term is the time derivative of an integral, which can be estimated against $\varepsilon^2 \mathcal{O}(\mathcal{E}_{\ell})$ by using Cauchy-Schwarz.
The last three integrals can be estimated against $\veps^{2}\,\mathcal{O}(\mathcal{E}_{\ell}+1)$ since $\Vert \partial_t f \Vert_{\infty}= \mathcal{O}(1)$
and due to \eqref{RES1}.
\\
For the second integral, applying  \eqref{komtaylor} gives us
\begin{align*}
&-\frac{1}{2}  \varepsilon^{2}  \int_{\mathbb{R}} \big[i \omega, f\big]  \frac{\gamma}{i \omega} \partial_x^{\ell} R_{j} \, \partial_x^{\ell} R_{-j} \, dx
\\[2mm]
&
\qquad\qquad\qquad
=
-\frac{1}{2}  \varepsilon^{2} 
\sum_{n=1}^{\lceil \deg^*(\gamma) \rceil-1} \frac{1}{n!}
 \int_{\mathbb{R}}  (-i)^n \omega^{(n)}  \frac{\gamma}{ \omega} \partial_x^{\ell} R_{j} \, \partial_x^{\ell} R_{-j} \, \partial_x^n f \,dx
 \\[2mm] 
 & \qquad\qquad\qquad \quad
 + \mathcal{O}(\varepsilon^2)\,
 \Vert \mathcal{R}(\frac{\gamma}{ \omega} \partial_x^{\ell} R_{j}, f) \Vert_{L^2}  \Vert \partial_x ^{\ell} R_{-j}\Vert_{L^2},
\end{align*}
where with \eqref{komrestneu} we can estimate $\Vert \mathcal{R}(\frac{\gamma}{ \omega} \partial_x^{\ell} R_{j}, f) \Vert_{L^2}\Vert \partial_x ^{\ell} R_{-j}\Vert_{L^2}=\mathcal{O}(\mathcal{E}_{\ell}+1)$.
\\
The integrals in the third and the forth place
can be written as a sum of some $\veps^{3}\,\mathcal{O}(\mathcal{E}_{\ell}+1)$-terms and $m$ many integrals of the form
\begin{align*}
\varepsilon^2 \,\int_{\mathbb{R}}  \gamma_k \partial_x^{\ell} R_{p_k} \, \partial_x^{\ell} R_{q_k} \, f_k \,dx
\end{align*}
with $m$, $\gamma_k$, $f_k$, $p_k$ and $q_k$ just as in the statement.
One sees this by exploiting \eqref{int und theta-1}, Leibniz's rule, \eqref{LeibnizR}, \eqref{komtaylor} and \eqref{int und theta}. 
Since we have by assumption
\begin{align}
\label{iterationsPSI}
\Vert   \tilde{R}_{\psi} \Vert_{H^{ \lceil \deg^*(\rho) \rceil+\lceil\deg(\omega) \rceil +1}}+ \Vert  \partial_t \tilde{R}_{\psi} \Vert_{H^{\lceil \deg^*(\rho) \rceil +1}}
&\le \varepsilon^{-1/2} c_{R} \,,
\\ \nonumber
\Vert  \tilde{R}_{\psi} \Vert_{C^{\lceil \deg^*(\rho) \rceil +\lceil\deg(\omega) \rceil}}
+\Vert  \partial_t  \tilde{R}_{\psi} \Vert_{C^{\lceil \deg^*(\rho) \rceil }} &\le c_{R} \,,
\end{align}
for some $c_{R} \in \mathbb{R}$,
straightforward estimates confirm that 
the
functions having the form $f_k= \varepsilon \,\partial_x^p \tilde{R}_{\psi}  f$ with $p\ge 0$  indeed fulfill 
\eqref{umform L^2}, \eqref{umform L^infty} and hence \eqref{gutes fk}.

%\begin{align*}
%&\Vert  M_{\varphi}^1f \Vert_{H^{n}}+ \Vert   \partial_t  M_{\varphi}^1 f \Vert_{H^m}
%\\
%& \qquad \qquad \qquad 
%\le   C \, \big( \Vert  f \Vert_{H^n} +\Vert \partial_t f \Vert_{H^m} \big)
%%+	   C \, \big( \Vert  \partial_x^{-1}f \Vert_{C^m} +\Vert  \partial_t \partial_x^{-1}f \Vert_{C^m} \big),
%%\\ \textit{and}
%\\
%& \Vert  \partial_x^{-1} M_{\varphi}^1 f \Vert_{\infty}+\Vert   \partial_t \partial_x^{-1} M_{\varphi}^1f\Vert_{C^{m}}
%\\
%& \qquad \qquad \qquad 
% \le  
% \varepsilon^{-1/2} \,
%  C \, \big( \Vert  f \Vert_{L^2}+\Vert \partial_t f \Vert_{L^2}  \big)
% +
% C \,   \big( \Vert  f \Vert_{C^{m-1}} +\Vert  \partial_t f \Vert_{C^{m-1}}  \big)
%,
%  \end{align*}
% for $m,n$ as required.

%and 
%\begin{align*}
%&\Vert M_{\varphi}^2 f \Vert_{H^{n}}+ \Vert  \partial_t  M_{\varphi}^2 f \Vert_{H^m}
%\\
%& \qquad \qquad \qquad 
%\le  C \, \big( \Vert  f \Vert_{H^n} +\Vert \partial_t f \Vert_{H^m} \big)+
%	 \varepsilon^{-1/2} \, C \, \big( \Vert  \partial_x^{-1}f \Vert_{C^m} +\Vert  \partial_t \partial_x^{-1}f \Vert_{C^m} \big),
%%\\ \textit{and}
%\\
%& \Vert  \partial_x^{-1} M_{\varphi}^2 f \Vert_{\infty}+\Vert  \partial_t \partial_x^{-1} M_{\varphi}^2 f\Vert_{C^{m}}
%\\
%& \qquad \qquad \qquad 
% \le  
%   C \,   \big( \Vert  \partial_x^{-1}f \Vert_{C^m} +\Vert  \partial_t \partial_x^{-1}f \Vert_{C^m}  \big)
%  .
%\end{align*}
% for
%\begin{align}
%\label{map2}
%&M^2_{\varphi}:f \mapsto \partial_x(\varphi \, \partial_x^{-1} f) \,.
%\end{align}
%This gives us that this 

\medskip
\textbf{b)}
Handling integrals of the form 
\begin{align}
\label{int rjrj skew}
\varepsilon^{2}  \int_{\mathbb{R}}& i \sigma \partial_x^{\ell} R_{j} \, \partial_x^{\ell} R_{j} \,   h \,dx \,.
\end{align}
%with \eqref{h prop} and \eqref{f prop}.
%\\
Since $i \sigma$ is skew symmetric and due to \eqref{komtaylor} and \eqref{komrestneu} we have
\begin{align*}
\varepsilon^{2}  \int_{\mathbb{R}} i \sigma \partial_x^{\ell} R_{j} \, \partial_x^{\ell} R_{j} \,   h \,dx
&=- \frac{1}{2}
\varepsilon^{2}  \int_{\mathbb{R}} \big[i \sigma, h \big] \partial_x^{\ell} R_{j} \, \partial_x^{\ell} R_{j} \,  dx,
\\[2mm]
&=
\varepsilon^2 \sum_{k=1}^{\deg^*(\sigma)-1}
\int_{\mathbb{R}}  {\varsigma}_k \partial_x^{\ell} R_{j} \, \partial_x^{\ell} R_{j} \, \partial_x^k h \,dx
+\veps^{2}\,\mathcal{O}(\mathcal{E}_{\ell}+1) \,,
\end{align*}
%what is the sum of integrals of the form  \eqref{int rjrj sym} and \eqref{int rjrj skew}, which have at least a whole derivative less than the original term,
%plus some $\veps^{2}\,\mathcal{O}(\mathcal{E}_{\ell}+1)$-terms 
with ${\varsigma}_k$ just as in the statement.

\medskip
\textbf{c)}
Handling integrals of the form
\begin{align}
\label{int rjrj sym}
\varepsilon^{2}  \int_{\mathbb{R}}& \upsilon \partial_x^{\ell} R_{j} \, \partial_x^{\ell} R_{j} \,   g \,dx \,.
\end{align}
%with \eqref{g prop} and \eqref{f prop}.
%\\
By using \eqref{komtaylor}, we can write
\begin{align*}
\varepsilon^{2}  \int_{\mathbb{R}}& \upsilon \partial_x^{\ell} R_{j} \, \partial_x^{\ell} R_{j} \,   g \,dx
\\[2mm]
& =
\varepsilon^{2}  \int_{\mathbb{R}} \big[i \omega, \partial_x^{-1} g\big]  \frac{\upsilon}{\omega'} \partial_x^{\ell} R_{j} \, \partial_x^{\ell} R_{j} \, dx
\\[2mm]
& \quad
+
\varepsilon^{2} 
\, 
\sum_{n=2}^{\lceil \deg^*(\upsilon)\rceil}
\frac{(-1)^{n}}{(n)!} \int_{\mathbb{R}} i^{n+1} \omega^{(n)}  \frac{\upsilon}{\omega'} \partial_x^{\ell} R_{j} \, \partial_x^{\ell} R_{j} \,\partial_x^{n-1} g \, dx
\\[2mm]
& \quad
 + \mathcal{O}(\varepsilon^2) \,
 \Vert \mathcal{R}(\frac{\upsilon}{ \omega'} \partial_x^{\ell} R_{j}, \partial_x^{-1} g) \Vert_{L^2}  \Vert \partial_x ^{\ell} R_{j}\Vert_{L^2}
 \, ,
%
%  alte Version
%-
%\varepsilon^{2}  \int_{\mathbb{R}} \mathcal{R} \big( \partial_x^{-1} g, \frac{\upsilon}{\omega'} \partial_x^{\ell} R_{j} \big) \, \partial_x^{\ell} R_{j} \, dx
%,
%\\
%& \quad
\end{align*}
where with \eqref{komrestneu} we can estimate  $\Vert \mathcal{R}(\frac{\upsilon}{ \omega'} \partial_x^{\ell} R_{j}, \partial_x^{-1} g) \Vert_{L^2} \Vert \partial_x ^{\ell} R_{j}\Vert_{L^2}=
\mathcal{O}(\mathcal{E}_{\ell}+1)$.
\\
Now, the second term already has the desired form and the last term is  $\varepsilon^{2}\,\mathcal{O}(\mathcal{E}_{\ell}+1)$ such that we only have to look at the first term.
%Since the second term is a sum of integrals of the form \eqref{int rjrj skew} and \eqref{int rjrj sym} which have at least
%a whole derivative less than our original term and the last term is  $\varepsilon^{2}\,\mathcal{O}(\mathcal{E}_{\ell}+1)$, we only have to look at the first term.
\\
By  exploiting the skew symmetry of $i \omega$ and \eqref{R_j} (and \eqref{RES1}), we have
\begin{align*}
&\varepsilon^{2}  \int_{\mathbb{R}} \big[i \omega, \partial_x^{-1}g\big]  \frac{\upsilon}{\omega'} \partial_x^{\ell} R_{j} \, \partial_x^{\ell} R_{j} \, dx
\\[2mm]
\nonumber
&
\qquad\qquad
=
\varepsilon^{2}  \int_{\mathbb{R}} i \omega \big( \partial_x^{-1}g  \frac{\upsilon}{\omega'} \partial_x^{\ell} R_{j} \big)\, \partial_x^{\ell} R_{j} \, dx
%\\[2mm] \nonumber
%& 
%\qquad\qquad\quad
-
\varepsilon^{2}  \int_{\mathbb{R}}  \partial_x^{-1}g \, \, i \omega  \frac{\upsilon}{\omega'} \partial_x^{\ell} R_{j} \, \partial_x^{\ell} R_{j} \, dx
\\[2mm]
\nonumber
&
\qquad\qquad
=-
j\,\varepsilon^{2} \, \partial_t \int_{\mathbb{R}}   \frac{\upsilon}{\omega'} \partial_x^{\ell} R_{j} \, \partial_x^{\ell} R_{j} \, \partial_x^{-1}g \,dx
\\[2mm] \nonumber
& \qquad\qquad\quad
+
\varepsilon^{3} \,  \int_{\mathbb{R}}   \frac{\upsilon}{\omega'} i \rho \partial_x^{\ell} \vartheta^{-1} \big(R_{\psi}  \vartheta (R_{-1}+R_{1}) \big) \, \partial_x^{\ell} R_{j} \,
\partial_x^{-1}g \, dx
\\[2mm] \nonumber
& \qquad\qquad\quad
+
\varepsilon^{3} \,  \int_{\mathbb{R}}   \frac{\upsilon}{\omega'} \partial_x^{\ell} R_{j} \, i \rho \partial_x^{\ell} \vartheta^{-1} \big(R_{\psi}  \vartheta (R_{-1}+R_{1}) \big) \, \partial_x^{-1}g \,dx
\\[2mm] \nonumber
& \qquad\qquad\quad
+
j\,\varepsilon^{2} \,  \int_{\mathbb{R}}   \frac{\upsilon}{\omega'} \partial_x^{\ell} R_{j} \, \partial_x^{\ell} R_{j} \, \partial_t \partial_x^{-1}g \,dx
\\[2mm] \nonumber
& \qquad\qquad\quad
+
\varepsilon^{2}\,\mathcal{O}(\mathcal{E}_{\ell}+1) \, .
\end{align*}
The last integral can be estimated against $\varepsilon^{2}\,\mathcal{O}(\mathcal{E}_{\ell}+1)$
since $\Vert \partial_t \partial_x^{-1}g \Vert_{\infty}= \mathcal{O}(1)$ .
Due to \eqref{int und theta-1}, the skew symmetry of $i \rho$ and the symmetry of $\omega'$ and $\upsilon$, we get
\begin{align*}
\varepsilon^{2}  \int_{\mathbb{R}}& \big[i \omega, \partial_x^{-1} g\big]  \frac{\upsilon}{\omega'} \partial_x^{\ell} R_{j} \, \partial_x^{\ell} R_{j} \, dx
\\[2mm]
\nonumber
&=-
j\,\varepsilon^{2} \, \partial_t \int_{\mathbb{R}}   \frac{\upsilon}{\omega'} \partial_x^{\ell} R_{j} \, \partial_x^{\ell} R_{j} \, \partial_x^{-1}g \,dx
\\[2mm] \nonumber
& \quad
-2
\varepsilon^{3} \,  \int_{\mathbb{R}}    i \rho \frac{\upsilon}{\omega'} \partial_x^{\ell} R_{j} \, \partial_x^{\ell} \big(R_{\psi} \vartheta(R_{-1}+R_{1}) \big) \, 
\partial_x^{-1}g \, dx
\\[2mm] \nonumber
& \quad
-
\varepsilon^{3} \,  \int_{\mathbb{R}}  \big[  i \rho \frac{\upsilon}{\omega'}, \partial_x^{-1}g \big] \partial_x^{\ell} R_{j} \,
 \partial_x^{\ell} \big(R_{\psi} \vartheta(R_{-1}+R_{1}) \big) \,  dx
\\[2mm] \nonumber
& \quad
-
\varepsilon^{3} \,  \int_{\mathbb{R}}  \big[i \rho,\partial_x^{-1}g  \big] \frac{\upsilon}{\omega'} \partial_x^{\ell} R_{j} \,  \partial_x^{\ell} \big(R_{\psi} \vartheta(R_{-1}+R_{1}) \big) \, dx
%\\[2mm] \nonumber
%& \quad
%+
%j\,\varepsilon^{2} \,  \int_{\mathbb{R}}   \frac{\upsilon}{\omega'} \partial_x^{\ell} R_{j} \, \partial_x^{\ell} R_{j} \, \partial_t \partial_x^{-1}g \,dx
\\[2mm] \nonumber
& \quad
+
\varepsilon^{2}\,\mathcal{O}(\mathcal{E}_{\ell}+1)
.
\end{align*}
The first term is a time derivative of an expression $\varepsilon^2 \tilde{D}$, which can be estimated against $\veps\,\mathcal{O}(\mathcal{E}_{\ell})$  since $ \varepsilon \Vert \partial_x^{-1}g\Vert_{\infty}= \mathcal{O}(1)$.
By using \eqref{komtaylor} and Leibniz's rule, we can write the third and the fourth integral as a sum of 
some $\veps^{3}\,\mathcal{O}(\mathcal{E}_{\ell}+1)$-terms and integrals of the form
 \begin{align*}
\varepsilon^2 \,\int_{\mathbb{R}}  \gamma_k \partial_x^{\ell} R_{p_k} \, \partial_x^{\ell} R_{q_k} \, f_k \,dx
\end{align*} 
with $\gamma_k$, $f_k$, $p_k$, $q_k$  just as in the statement. Making straightforward estimates by using \eqref{iterationsPSI}
shows that the functions of the form $f_k= \varepsilon \,\partial_x^n \tilde{R}_{\psi} \partial_x^m g$ with $n,m \ge 0$
 here fulfill \eqref{gutes fk}, \eqref{umform L^2} and \eqref{umform L^infty}. 
\\
What remains to be analyzed is the second term.
Using Leibniz's rule and afterwards \eqref{LeibnizR} and \eqref{int und theta}, we obtain
\begin{align*}
-2
\varepsilon^{3} \,  \int_{\mathbb{R}}&    i \rho \frac{\upsilon}{\omega'} \partial_x^{\ell} R_{j} \, \partial_x^{\ell} \big(R_{\psi} \vartheta (R_{-1}+R_{1}) \big) \, 
\partial_x^{-1}g \, dx
\\[2mm] \nonumber
&=
-2
\varepsilon^{3} \,  \int_{\mathbb{R}}    i \rho \frac{\upsilon}{\omega'} \partial_x^{\ell} R_{j} \, \partial_x^{\ell}(R_{-1}+R_{1}) \, \tilde{R}_{\psi} \, 
\partial_x^{-1}g \, dx
\\[2mm] \nonumber
& \quad
-2\varepsilon^{3} \,
\sum_{m=1}^{M} 
\binom{\ell}{m} \,
  \int_{\mathbb{R}}    i \rho \frac{\upsilon}{\omega'} \partial_x^{\ell} R_{j} \, \partial_x^{\ell-m}  (R_{-1}+R_{1})  \, \partial_x^{m} \tilde{R}_{\psi} \,
\partial_x^{-1}g \, dx
\\[2mm] \nonumber
& \qquad
+
\veps^{3}\,\mathcal{O}(\mathcal{E}_{\ell}+1),
\end{align*}
where $M:=  \lceil \deg^*(\rho \upsilon)- \deg (\omega')\rceil -1$.
\\
The second term here
consists (after integration by parts) of integrals having the form
 \begin{align*}
\varepsilon^2 \,\int_{\mathbb{R}}  \gamma_k \partial_x^{\ell} R_{p_k} \, \partial_x^{\ell} R_{q_k} \, f_k \,dx
\end{align*} 
with $\gamma_k$, $f_k$, $p_k$, $q_k$  just as in the statement. The functions of the form $f_k= \varepsilon \,\partial_x^m \tilde{R}_{\psi} \partial_x^{-1} g$ with $m > 0$ 
 fulfill  \eqref{umform L^2} and \eqref{umform L^infty}:
For $m \ge 1$ and $n$ as required one gets
\begin{align*}
\Vert \varepsilon \partial_x^m \tilde{R}_{\Psi} \partial_x^{-1}g\Vert_{H^{n}}
\le\varepsilon \, \Vert  \partial_x^m \tilde{R}_{\Psi} \Vert_{H^{n}} \Vert \partial_x^{-1}g\Vert_{C^{n}}
\,,
\\
\Vert \varepsilon \partial_x^m \tilde{R}_{\Psi} \partial_x^{-1}g\Vert_{C^{n}}
\le\varepsilon \, \Vert  \partial_x^m \tilde{R}_{\Psi} \Vert_{C^{n}} \Vert \partial_x^{-1}g\Vert_{C^{n}} \,.
\end{align*}
The estimates  for $\Vert \varepsilon  \partial_t (\partial_x^m \tilde{R}_{\Psi} \partial_x^{-1}g )\Vert_{H^{n}}$
and $\Vert \varepsilon  \partial_t (\partial_x^m \tilde{R}_{\Psi} \partial_x^{-1}g )\Vert_{C^{n}}$
are similarly straightforward. Concerning the other estimates, we estimate 
\begin{align*}
\Vert \varepsilon  \partial_x^{-1} (\partial_x^m \tilde{R}_{\Psi} \partial_x^{-1}g)\Vert_{\infty}
&= \varepsilon 
\Vert \partial_x^{m-1} \tilde{R}_{\Psi} \partial_x^{-1}g
-\partial_x^{-1} (\partial_x^{m-1} \tilde{R}_{\Psi} g)
\Vert_{\infty}
\\
&\le 
\varepsilon
\Vert \partial_x^{m-1} \tilde{R}_{\Psi}\Vert_{\infty} \Vert \partial_x^{-1}g\Vert_{\infty}
+ \varepsilon 
\Vert\partial_x^{m-1} \tilde{R}_{\Psi} \Vert_{L^2} \Vert g \Vert_{L^2} \,,
\end{align*}
and similar $\Vert \varepsilon  \partial_t \partial_x^{-1} (\partial_x^m \tilde{R}_{\Psi} \partial_x^{-1}g)\Vert_{\infty}$.
The  estimate
$\Vert \varepsilon \partial_x^m \tilde{R}_{\Psi} \partial_x^{-1}g\Vert_{H^{n}}= \mathcal{O}(1)$ 
is not implied by the above estimates
such that \eqref{gutes fk} is not trivially obtained. 
We confirm  \eqref{gutes fk}  by exploiting that the supremum over all $x \in \mathbb{R}$ is the same as the supremum over all 
$ (\varepsilon^{-1} x) \in \mathbb{R}$ such that the loss of $\varepsilon$-powers caused by the slow spatial scale of the NLS  present in  estimate \eqref{iterationsPSI} can here be avoided.
%This is since $\tilde{R}_{\Psi}= \psi_c+ \varepsilon \tilde{R}_{Q}$ and $\Vert \psi_c \Vert_{L^2}= \mathcal{O}(\varepsilon^{-1/2}) $ due the  scaling
%of the NLS approximation. We still can explicitly calculate that \eqref{gutes fk} holds by exploiting that the supremum over all $x \in \mathbb{R}$ is the same as the supremum over all $ (\varepsilon^{-1} x) \in \mathbb{R}$.
\\
Thus, we now only have to examine the term
\begin{align*}
-2
\varepsilon^{3} \,  \int_{\mathbb{R}}&    i \rho \frac{\upsilon}{\omega'} \partial_x^{\ell} R_{j} \, \partial_x^{\ell}(R_{-1}+R_{1}) \, \tilde{R}_{\psi} \, 
\partial_x^{-1}g \, dx
\\[2mm] \nonumber
&=
-2\varepsilon^{3} \,  \int_{\mathbb{R}}    i \rho \frac{\upsilon}{\omega'} \partial_x^{\ell} R_{j} \, \partial_x^{\ell} R_{-j} \, \tilde{R}_{\psi} \, \partial_x^{-1}g \, dx
\\
& \quad
-2\varepsilon^{3} \,  \int_{\mathbb{R}}    i \rho \frac{\upsilon}{\omega'} \partial_x^{\ell} R_{j} \, \partial_x^{\ell} R_{j} \, \tilde{R}_{\psi} \, \partial_x^{-1}g \, dx
\,.
\end{align*}
For the above first integral
 we can proceed as in  paragraph \textbf{a)}
 and for the second integral as in  paragraph \textbf{b)}.
The required estimates for the function $\tilde{h}=\varepsilon \tilde{R}_{\psi} \, \partial_x^{-1}g$ work as above, however one has to be aware  
that  $\tilde{h}$ only meets the
 conditions for \textbf{a)} and \textbf{b)} since
$
\deg^*(\rho \frac{\upsilon}{\omega'}) \le \deg^*(\upsilon)+1 
$ and $
\deg^*(\rho \frac{\upsilon}{\omega'}) \le \deg^*(\rho) \le \deg^*(\omega)$ due to \eqref{rho BED2}. Also note that for \textbf{a)} and \textbf{b)} no estimate for $\partial_x^{-1}\tilde{h}$ is needed.
\\
With this the statement is proven.

\medskip
One important aspect of the now valid  statement is
that we can use it repeatedly, i.e. we can always apply it again for every integral on the right hand side
of \eqref{umform}.
So we can use \eqref{umform} and exploit \eqref{umform ord'} in order to get
\begin{align*}
\varepsilon^{2}  \int_{\mathbb{R}}&  \gamma \partial_x^{\ell} R_{j_1} \, \partial_x^{\ell} R_{j_2} \, f \,dx
\\ \nonumber
&
=\varepsilon^2 \partial_t \tilde{\mathcal{D}}
+\varepsilon^2 \sum_{k=1}^{\tilde{m}_{\gamma}}
\int_{\mathbb{R}}  \tilde{\gamma}_k \partial_x^{\ell} R_{p_k} \, \partial_x^{\ell} R_{q_k} \, \tilde{f}_k \,dx
+\varepsilon^{2}\,\mathcal{O}(\mathcal{E}_{\ell}+1),
\end{align*}
where $\tilde{\mathcal{D}}$ with $\varepsilon^{2} \,\tilde{\mathcal{D}}= \varepsilon \, \mathcal{O}(\mathcal{E}_{\ell})$, 
$\tilde{m}_{\gamma} \le m= m(\deg(\omega)) \in \mathbb{N}$
and due to \eqref{umform L^2} and \eqref{umform L^infty}:
\begin{align*}
\Vert \tilde{f}_k\Vert_{H^{\lceil \deg^*(\tilde{\gamma}_k) \rceil \lceil \deg(\omega) \rceil }}+ \Vert \partial_t \tilde{f}_k \Vert_{H^{\lceil \deg^*(\tilde{\gamma}_k) \rceil}}
&\le \varepsilon^{-1/2}\, C_{f},
\\
\Vert \partial_x^{-1} \tilde{f}_k\Vert_{\infty} +\Vert \partial_t \partial_x^{-1} \tilde{f}_k\Vert_{C^{\lceil \deg^*(\tilde{\gamma}_k) \rceil}}
&\le C_{f} \,,
\end{align*}
for some constant $C_{f}=C_{f}(\tilde{R}_{\psi}, f, \gamma)> 1$.
\\
By using \eqref{umform} and exploiting  \eqref{umform ord'} again for every integral on the above right-hand side, we can obtain
an expression $\varepsilon^{2} \, \check{\mathcal{D}}=\varepsilon \,\mathcal{O}(\mathcal{E}_{\ell})$ 
such that we have
\begin{align*}
\varepsilon^{2}  \int_{\mathbb{R}}&  \gamma \partial_x^{\ell} R_{j_1} \, \partial_x^{\ell} R_{j_2} \, f \,dx
\\ \nonumber
&
=\varepsilon^2 \partial_t \check{\mathcal{D}}
+\varepsilon^2 \sum_{k=1}^{\check{m}}
\int_{\mathbb{R}}  \check{\gamma}_k \partial_x^{\ell} R_{p_k} \, \partial_x^{\ell} R_{q_k} \, \check{f}_k \,dx
+\varepsilon^{2}\,\mathcal{O}(\mathcal{E}_{\ell}+1),
\end{align*}
where  $\check{m} \le m^2$, 
\begin{align*}
\Vert \check{f}_k\Vert_{H^{\lceil \deg^*(\check{\gamma}_k) \rceil \lceil \deg(\omega) \rceil }}+ \Vert \partial_t \check{f}_k \Vert_{H^{\lceil \deg^*(\check{\gamma}_k) \rceil}}
&\le
\varepsilon^{1/2} \, C_{f}^2 \,,
\\
\Vert \partial_x^{-1} \check{f}_k\Vert_{\infty} +\Vert \partial_t \partial_x^{-1} \check{f}_k\Vert_{C^{\lceil \deg^*(\check{\gamma}_k) \rceil}}
&\le C_{f}^2 \,.
\end{align*}
By repeating the last step  $N+1$ times, we get 
\begin{align*}
\varepsilon^{2}  \int_{\mathbb{R}}&  \gamma \partial_x^{\ell} R_{j_1} \, \partial_x^{\ell} R_{j_2} \, f \,dx
\\ \nonumber
&
=\varepsilon^2  \sum_{p=0}^N \varepsilon^{p/2} \, \partial_t{\mathcal{D}}_p
+\varepsilon^{2} \, \varepsilon^{\frac{N+1}{2}}   \sum_{k=1}^{m_N}
\int_{\mathbb{R}}  \gamma_{k,N} \partial_x^{\ell} R_{p_k} \, \partial_x^{\ell} R_{q_k} \, f_{k,N} \,dx
+\varepsilon^{2} \sum_{p=0}^N \varepsilon^{p/2} \,\mathcal{C}_p,
\end{align*}
for some expressions ${\mathcal{D}}_p$ with $\varepsilon^{2} \,{\mathcal{D}}_p=\varepsilon \,\mathcal{O}(\mathcal{E}_{\ell})$, some $\mathcal{C}_p=\mathcal{O}(\mathcal{E}_{\ell}+1)$, 
$m_N \le m^{3+N}$ and 
\begin{align*}
\Vert f_{k,N}\Vert_{H^{\lceil \deg^*(\gamma_{k,N}) \rceil+\lceil \deg(\omega) \rceil}}+ \Vert \partial_t f_{k,N} \Vert_{H^{\lceil \deg^*(\gamma_{k,N}) \rceil}}
&\le 
\varepsilon^{1/2} \, C_{f}^{3+N},
\\
\Vert \partial_x^{-1} f_{k,N} \Vert_{\infty} +\Vert \partial_t \partial_x^{-1} f_{k,N}\Vert_{C^{\lceil \deg^*(\gamma_{k,N}) \rceil}}
&\le C_{f}^{3+N} .
\end{align*}
Moreover, we have $\deg^*(\gamma_{k,N})\le \deg^*(\rho)$ .
%\textbf{Bemerkung}
%\\
%\textit{Wir haben dabei sogar $\deg^*(\gamma_{k,N}) \in \{\deg^*(\gamma), \deg^*(\gamma)-1,\deg^*(\gamma)-2, \dots \}$. }
%
\\
We will now  show  that 
\begin{align*}
\mathcal{D}^{\infty}:=
\sum_{p=0}^{\infty} \varepsilon^{p/2} \, {\mathcal{D}}_p
\end{align*}
does exist,  $\varepsilon^{2} \,\mathcal{D}^{\infty}=\varepsilon \,\mathcal{O}(\mathcal{E}_{\ell})$  and  
\begin{align*}
\varepsilon^{2}  \int_{\mathbb{R}}  \gamma \partial_x^{\ell} R_{j_1} \, \partial_x^{\ell} R_{j_2} \, f \,dx
=\varepsilon^2  \, \partial_t \mathcal{D}^{\infty}
+\varepsilon^{2}  \,\mathcal{O}(\mathcal{E}_{\ell}+1).
\end{align*}
%First off, we have 
%\begin{align*}
%&\varepsilon^{\frac{N+1}{2}}   \sum_{k=1}^{m_N}
%\int_{\mathbb{R}}  \gamma_{k,N} \partial_x^{\ell} R_{p_k} \, \partial_x^{\ell} R_{q_k} \, f_{k,N} \,dx
%\\[2mm]
%& \qquad\le \varepsilon^{\frac{N+2}{2}} \, C_{f}^{N+3}
%\Big(
%\Vert  R_{1} \Vert_{H^{\ell}}\Vert  R_{1} \Vert_{C^{\ell+ \lceil \deg^*(\gamma) \rceil}} 
%+
%\Vert  R_{1} \Vert_{H^{\ell}}\Vert  R_{-1} \Vert_{C^{\ell+ \lceil \deg^*(\gamma) \rceil}}
%\\
%& \qquad \qquad \qquad   \qquad \qquad \quad
%+
%\Vert  R_{-1} \Vert_{H^{\ell}}\Vert  R_{1} \Vert_{C^{\ell+ \lceil \deg^*(\gamma) \rceil}}
%+
%\Vert  R_{-1} \Vert_{H^{\ell}}\Vert  R_{-1} \Vert_{C^{\ell+ \lceil \deg^*(\gamma) \rceil}}
%\Big).
%\end{align*}
By taking a close look at the proof of \eqref{umform}, we find that
\begin{align*}
\varepsilon^{2} \,
\varepsilon^{\frac{p}{2}}\,{\mathcal{D}}_p  &\le 
\varepsilon \,\varepsilon^{\frac{p}{2}} \,  m^{p+3} c^{p+3}  C_{f}^{p+3} \, \mathcal{E}_{\ell},
\\
 \varepsilon^{\frac{p}{2}} \,{\mathcal{C}}_p &\le \varepsilon^{\frac{p}{2}} \, m^{p+3} c^{p+3}  C_{f}^{p+3} \, \big(\mathcal{E}_{\ell}+1 \big),
\end{align*}
for some  $c>1$  as long as   $f, i \rho, i \omega$ and $\ell$ are fixed.
We emphasize that this is in particular possible  due to the fact that $\deg^*(\gamma_{k,N})$ is always uniformly bounded by $\deg^*(\omega)$.
\\
By now choosing $\varepsilon_0$ small enough,
for instance such that
\begin{align*}
\varepsilon_{0}^{1/4}\,m c \, C_{f} \le 1 \,,
\end{align*}
we get the following.
There is a $c \in \mathbb{R}$ such that
\begin{align*}
\varepsilon^{2} \,\mathcal{D}^{\infty} 
= \varepsilon^{2} \,\sum_{p=0}^{\infty} \varepsilon^{p/2} \,  {\mathcal{D}}_p 
\le\varepsilon^{2} \, \sum_{p=0}^{\infty} \varepsilon^{p/2} \, \vert {\mathcal{D}}_p \vert
\le \varepsilon \,
\sum_{p=0}^{\infty} \varepsilon^{p/4} \, c \,\mathcal{E}_{\ell} =
\varepsilon \,c \, \mathcal{E}_{\ell}\,  \, \sum_{p=0}^{\infty} \varepsilon^{p/4} \,
 =\varepsilon \,\mathcal{O}(\mathcal{E}_{\ell}) \,,
\end{align*}
analogously we get
\begin{align*}
\sum_{p=0}^{\infty} \varepsilon^{p/2} \, C_p \le \sum_{p=0}^{\infty} \varepsilon^{p/2} \, \vert C_p \vert= \mathcal{O}(\mathcal{E}_{\ell}+1) \,.
\end{align*}
Moreover,
%, assuming $R_{-1}, R_1 \in  H^{\ell} \cap {C}_b^{\ell+ \lceil\deg^*(\rho)\rceil}$, 
\begin{align*}
&\varepsilon^{\frac{N+1}{2}}   \sum_{k=1}^{m_N}
\int_{\mathbb{R}} 
  \gamma_{k,N} \partial_x^{\ell} R_{j_k} \, \partial_x^{\ell} R_{l_k} \, f_{k,N} \,dx
\\[2mm]
& \qquad \qquad
  \le 
  \varepsilon^{\frac{N+1}{4}}  \, C_{f}^2 \,
\Big(
\Vert  R_{1} \Vert_{H^{\ell}}\Vert  R_{1} \Vert_{C^{\ell+ \lceil \deg^*(\gamma) \rceil}} 
+
\Vert  R_{1} \Vert_{H^{\ell}}\Vert  R_{-1} \Vert_{C^{\ell+ \lceil \deg^*(\gamma) \rceil}}
\\
& \qquad \qquad \qquad \qquad \qquad
+
\Vert  R_{-1} \Vert_{H^{\ell}}\Vert  R_{1} \Vert_{C^{\ell+ \lceil \deg^*(\gamma) \rceil}}
+
\Vert  R_{-1} \Vert_{H^{\ell}}\Vert  R_{-1} \Vert_{C^{\ell+ \lceil \deg^*(\gamma) \rceil}}
\Big)
\\[2mm]
& \qquad \qquad
=0, \quad \text{for  } N \rightarrow \infty.
\end{align*}
The short involvement of the $C^{\ell+\lceil\deg^*\rho \rceil}$-norm here is not problematic since the final estimate does no longer involve this norm.
\\
We now obtain
\begin{align*}
\varepsilon^{2}  \int_{\mathbb{R}}  \gamma \partial_x^{\ell} R_{j_1} \, \partial_x^{\ell} R_{j_2} \, f \,dx
&=\varepsilon^2  \, \partial_t \mathcal{D}^{\infty}
+\varepsilon^{2} \sum_{p=0}^{\infty} \varepsilon^{p/2} \, C_p
\\
& \qquad 
+ \varepsilon^2
\lim_{N \rightarrow \infty}
\varepsilon^{\frac{N+1}{2}}   \sum_{k=1}^{m_N}
\int_{\mathbb{R}} 
  \gamma_{k,N} \partial_x^{\ell} R_{q_k} \, \partial_x^{\ell} R_{p_k} \, f_{k,N} \,dx
\\[2mm]
&=\varepsilon^2  \, \partial_t \mathcal{D}^{\infty}
+\varepsilon^{2}  \,\mathcal{O}(\mathcal{E}_{\ell}+1).
\end{align*}
\qed
\\
For $\deg^* (\rho) \le 1$ or $\deg^* (\rho) < \deg(\omega)$, one can modify the above proof by exploiting the special structure of \eqref{QDS},
i.e. 
\begin{align}
\label{d/dt (r1+r-1)}
\partial_t(R_1 +R_{-1})
&=  i \omega (R_{1}- R_{-1})  +\varepsilon^{-\beta} \vartheta^{-1} \big(\mathrm{Res}_{u_1}(\varepsilon \Psi)+\mathrm{Res}_{u_{-1}}(\varepsilon \Psi) \big) \,,
\end{align}
to obtain a terminating algorithm that gives out an explicit expression $\varepsilon^2 \mathcal{D}$ consisting of a finite sum of integrals.

\medskip
Corollary \ref{corollary checkEnergy} now allows us to prove Theorem \ref{mainresult}.
\\
\textbf{Proof of Theorem \ref{mainresult} }
For $\ell \ge \lceil \deg(\omega) \rceil + \lceil \deg^*(\rho) \rceil +1$,
we can use corollary \ref{corollary checkEnergy} together with
Gronwall's inequality in order to obtain the $\mathcal{O}(1)$-boundedness of $\tilde{\mathcal{E}_{\ell}}$ for all $t\in[0,T_0/\veps^{2}]$
as long as  $\varepsilon_0>0$ is chosen sufficiently small.
%\\
%Since we assumed the local existence of solutions to \eqref{QDS}, there is a some $T(\varepsilon) >0$ such that the $H^{\ell}$-norms of $R_{-1}(t)$ and $R_{1}(t)$ can be  uniformly bounded as long as $0\le t \le T(\varepsilon) $. 
%Due to corollary \ref{cor32}, we know that for sufficiently small $\varepsilon_0$, \eqref{epsilon E_l le 1} is true for $0\le t \le T(\varepsilon) $, i.e.
%\begin{align*}
%\varepsilon \mathcal{E}_{\ell}(t) \le 1 \,,
%\end{align*}
% and thus corollary \ref{corollary checkEnergy} does indeed hold for 
%$0\le t \le T(\varepsilon) $.
%\\
%In particular, we have 
%\begin{align*}
%\partial_t \tilde{\mathcal{E}}_{\ell}(t) \le \varepsilon^2 \, C \, \big( \tilde{\mathcal{E}}_{\ell}(t)+1 \big) 
%\, 
%\end{align*}
%for some $C \ge 0$ and $0\le t \le T(\varepsilon) $.
%\\
%Gronwall's inequality now yields
%\begin{align*}
%\tilde{\mathcal{E}}_{\ell}(t) \le  \big( \tilde{\mathcal{E}}_{\ell}(0) +\varepsilon^2 \,C t \big) e^{ \varepsilon^2 \, C t}
%\end{align*}
%for $0\le t \le T(\varepsilon) $.
%\\
%Choosing $\varepsilon_0$ such small that 
%\begin{align*}
%\big( \tilde{\mathcal{E}}_{\ell}(0) +C T_0 \big) e^{  C T_0} \le \varepsilon_0^{-1} \,,
%\end{align*}
%we can now obtain $T(\varepsilon) \ge \varepsilon^{-2} \, T_0 $, i.e. in particular
%\begin{align*}
%\tilde{\mathcal{E}}_{\ell}(t) \le  \big( \tilde{\mathcal{E}}_{\ell}(0) +\varepsilon^2 \,C t \big) e^{\varepsilon^2 \,C t}
%\end{align*}
%for $0\le t \le \varepsilon^{-2} \, T_0 $.
For sufficiently small $\varepsilon_0>0$ there thus is some constant $C_R$ such that
\begin{align*}
\sup_{[0, T_0/\varepsilon^2]}  
\Big\Vert \left(
\begin{array}{c} { R_{-1}} \\  {R}_1
\end{array}
\right) \Big\Vert_{{H^{\ell}}} \le C_R \, ,
\end{align*}
due to corollary \ref{corollary checkEnergy}.
Choosing $\ell \ge s_A$,  estimate \eqref{RES2} now allows to conclude
\begin{align*}
&\sup_{[0, T_0/\varepsilon^2]} \Vert
u- \varepsilon \psi_{NLS}\Vert_{H^{s_A}}
%=
%\sup_{[0, T_0/\varepsilon^2]} \Vert
%u_{-1}+ u_{1}- \varepsilon \psi_{NLS}\Vert_{H^{s_A}}
\\[2mm]
&
\qquad \qquad
 \le
\sup_{[0, T_0/\varepsilon^2]} \Big\Vert \left(
\begin{array}{c} {u_{-1}} \\ {u}_1
\end{array}
\right)
- \varepsilon \left(
\begin{array}{c}  \psi_{NLS} \\ 0
\end{array}
\right)  \Big\Vert_{{H^{s_A}}}
\\[2mm]
& \qquad \qquad 
\le\sup_{[0, T_0/\varepsilon^2]}  \varepsilon^{\beta} 
\Big\Vert \left(
\begin{array}{c} {\vartheta R_{-1}} \\ \vartheta {R}_1
\end{array}
\right) \Big\Vert_{{H^{s_A}}}
 +
\sup_{[0, T_0/\varepsilon^2]} 
 \Vert\varepsilon \Psi -\varepsilon \left(
\begin{array}{c}  \psi_{NLS} \\ 0
\end{array}
\right)  
 \Big\Vert_{{H^{s_A}}}
\\[2mm]
& \qquad \qquad \le \mathcal{O}( \varepsilon^{3/2}). 
\end{align*}	
\qed

\section{Discussion}
As model problem for the 2D water wave problem with finite depth and   surface tension $b \ge 0$,
one can look at \eqref{QDS} with 
\begin{align*}
\omega(k)&= \rho(k)=\sign(k) \sqrt{ k \tanh(k) (1+b k^2 ) } \,,
\end{align*}
%and  {$\rho(k)= \omega(k)$},
%or $\rho(k)= \sign(k) \sqrt{k  \tanh(k)} + b k$,
 cf. \cite{SW11, CW17} for the case without surface tension $b=0$.
Theorem \ref{mainresult} grants us the validity of the NLS approximation
for all $b \ge 0$ and $k_0 >0$, excluding some special pairs $(b,k_0)$ with $0<b<1/3$.
Indeed, the validity of the NLS approximation for the full  2D water wave problem with finite depth and surface tension was recently proven in \cite{D19}.
\\
System
\eqref{QDS}
with $\omega$ and $\rho$  given in Fourier space by 
\begin{align*}
\omega(k)= \rho(k)= \sign(k) \sqrt{ \frac{k \tanh(k)}{1 +k \tanh(k)} (1+k^2+k^4 ) } \,\,,
\end{align*}
can be considered as a model problem for the 2D water wave problem with finite depth  and ice cover.
This model has the same linear dispersion relation as the full problem, see e.g. \cite{I15}.
Moreover, its quasilinear quadratic term shares principle difficulties with the ones of the full problem
regarding the construction of the normal form transformations and the loss of regularity in the error estimates. 
%One could also include a dependence on the surface tension and the properties of the ice cover as in the full problem.
We omit an analysis of the possible resonances that can occur but
for some  $k_0$, e.g. $k_0=1$, the resonance condition of this paper is fulfilled such that  Theorem \ref{mainresult}  grants the validity of the NLS approximation for these wavenumbers.
Thus, the techniques of this paper might be useful 
%to handle the problematic quasilinear terms losing two  derivatives in 
for proving a NLS validity result for 
the full 2D water wave problem with ice cover.

Our result  could also be interesting for double dispersion equations.
With \eqref{beam equation}, we already gave one example but Theorem \ref{mainresult} also applies to other quasilinear double dispersion equations, like e.g.
$\partial_t^2 u = \partial_x^6 u +\partial_x^6 u^2 $ or
$\partial_t^2 u = -\partial_x^4 u + \partial_x^2 u + \partial_x^2 u^2  -\partial_x^4 u^2 $.

%Due to analogies in the methods of proof, it is possible to adjust the techniques of this paper  to
% prove existence of small, smooth solutions over cubically nonlinear time scales for quasilinear quadratic systems, cf. \cite{DH18}.

\medskip

\textbf{Acknowledgment:} This work was supported by the Deutsche Forschungsgemeinschaft DFG
under the grant DU 1198/2. The author thanks Guido Schneider and Wolf-Patrick D{\"u}ll for discussions.
\bibliographystyle{alpha}

\begin{thebibliography}{99}
%
%\bibitem[A01]{A01}
%Agrawal, G. P.:
%Solitons and nonlinear wave equations.
%Academic Press (2001).

%\bibitem[A03]{A03}
%Ambrose, D. M.:
%Well-posedness of vortex sheets with surface tension. 
%SIAM J. Math. Anal. {\bf 35}(1), 211-244 (2003). 



\bibitem[AS81]{AS81} Ablowitz, M.J., Segur, H.:  Solitons and the inverse scattering transform. In: SIAM Studies in Applied Mathematics, vol. 4. SIAM
(1981).

%\bibitem[C17]{C17}
%Cummings, P.:
%Nonlinear Schr{\"o}dinger approximations for partial differential equations with quadratic and quasilinear terms.
%PhD thesis,
%OpenBU: \url{https://open.bu.edu/handle/2144/24090} (2017).

%\bibitem[C87]{C87}
%Craig, W.: Nonstrictly hyperbolic nonlinear systems. {Math. Ann.}, {\bf 277}(2), 213-232 (1987).
%



%\bibitem[CDS15]{CDS15}
%Chirilus-Bruckner, M., D\"ull, W.-P., Schneider, G.:
%NLS approximation of time oscillatory long waves for equations with quasilinear quadratic terms.
%{Math. Nachr.} {\bf 288}(2-3), 158-166 (2015).

%\bibitem[CS13]{CS13}
%Chong, C., Schneider, G.:
%Numerical evidence for the validity of the NLS approximation in systems with a quasilinear quadratic nonlinearity.
%{ZAMM Z. Angew. Math. Mech.} {\bf 93}(9), 688-696 (2013).



\bibitem[CW17]{CW17}
Cummings, P., Wayne, C. E.:
Modified energy functionals and the NLS approximation.
{Discrete Contin. Dyn. Syst.} {\bf 37}(3), 1295-1321 (2017).

%\bibitem[D05]{D05}
%Debnath, L.:
%Nonlinear partial differential equations for scientists and engineers.
%Second ed., Birkhauser Boston Inc., Boston, MA, (2005).


%\bibitem[D12]{D12}
%D{\"u}ll, W.-P.:
%Validity of the Korteweg-de Vries Approximation for the Two-Dimensional Water Wave Problem in the Arc Length Formulation.
%{Comm. Pure Appl. Math.} {\bf 65}(3), 381-429 (2012).

\bibitem[D17]{D17}
D{\"u}ll, W.-P.:
Justification of the nonlinear Schr\"odinger approximation for a quasilinear Klein-Gordon equation.
Comm. Math. Phys. {\bf 355}(3), 1189-1207  (2017).
%
%\bibitem[D18]{D18}
%D{\"u}ll,  W.-P.: On the Mathematical Description of Time-Dependent Surface Water Waves.
%Jahresber. Dtsch. Math. Ver. {\bf 120} (2), 117-141 (2018).

\bibitem[D19]{D19}
D{\"u}ll,  W.-P.: Validity of the nonlinear Schr\"odinger approximation for the two-dimensional water wave problem with and without surface tension in the arc length
formulation.
 Arch. Ration. Mech. Anal.,  {https://doi.org/10.1007/s00205-020-01586-4} (2020).

\bibitem[DH18]{DH18}
D{\"u}ll, W.-P.; He{\ss}, M.:
Existence of long time solutions and validity of the nonlinear Schr\"odinger approximation for a quasilinear dispersive equation. 
J. Differential Equations {\bf 264}(4), 2598-2632 (2018). 

%\bibitem[DHSZ16]{DHSZ16}
%D{\"u}ll, W.-P., Hermann, A., Schneider, G., Zimmermann, D.:
%Justification of the 2D NLS equation -- Quadratic resonances do not matter 
%in case of analytic initial conditions. { J. Math. Anal. Appl.} {\bf 436}, 847-867 (2016).

\bibitem[DS06]{DS06}
D{\"u}ll, W.-P., Schneider, G.:
Justification of the Nonlinear Schr\"odinger equation 
for a resonant Boussinesq model.
{Indiana Univ. Math. J.} {\bf 55}(6), 1813-1834 (2006).

\bibitem[DSW16]{DSW16}
D{\"u}ll, W.-P., Schneider, G., Wayne, C.E.:
Justification of the Nonlinear Schr\"odinger equation for the evolution of gravity driven 2D surface water waves in a canal of finite depth. 
{Arch. Rat. Mech. Anal.} {\bf 220}(2), 543-602 (2016).

\bibitem[H19]{H19}
He{\ss}, M.:
Validity of the nonlinear Schr{\"o}dinger approximation for quasilinear dispersive systems,
PhD thesis, Universit{\"a}t Stuttgart (2019).

%\bibitem[HIT16]{HIT16}
%Hunter, J.K., Ifrim, M., Tataru, D.: 
%Two dimensional water waves in holomorphic coordinates,  Comm. Math. Phys. {\bf 346}(2),  483-552 (2016).

\bibitem[HITW15]{HITW15} J.K.Hunter, M.Ifrim, D.Tataru, T.K.Wong, Long Time Solutions for a Burgers-Hilbert Equation via a Modified Energy Method,
Proc. Amer. Math. Soc. {\bf 143}(8),  3407-3412 (2015). 

%\bibitem[IT15]{IT15}
%{Ifrim, M., Tataru, D.:}
%The lifespan of small data solutions in two dimensional capillary water waves, arXiv:1406.5471v2 (2014).


\bibitem[I15]{I15}
{Il'ichev, A. T.:}
Envelope solitary waves and dark solitons at a water-ice interface.
Proc. Steklov Inst. Math. {\bf 289}(1) , 152-166 (2015). 




\bibitem[IT19]{IT19}
{Ifrim, M., Tataru, D.:}
The NLS approximation for two dimensional deep gravity waves.
Sci. China Math. {\bf 62}(6), 1101-1120 (2019). 

\bibitem[K88]{K88}
Kalyakin, L.A.:
Asymptotic decay of a one-dimensional wave packet in a nonlinear
  dispersive medium.
{Sb. Math. } {\bf 60}, 457-483 (1988).

%\bibitem[K75]{K75}
%Kato, T.: The Cauchy problem for quasi-linear symmetric hyperbolic systems,
%{Arch. Rat. Mech. Anal.} {58},  181-205 (1975).

%\bibitem[K75b]{K75b}
%Kato, T.:Quasi-linear equations of evolution, with applications to partial differential equations.
%In: Everitt W.N. (eds) Spectral Theory and Differential Equations. Lecture Notes in Mathematics {\bf 448}, Springer, Berlin, Heidelberg,
%25-70 (1975). 


\bibitem[KV19]{KV19}
Kolkovska, N., Vucheva, V.: Invariant preserving schemes for double dispersion equations.
Adv. Difference Equ. {\bf 216} 1-16 (2019). 




%\bibitem[KSM92]{KSM92}
%Kirrmann, P., Schneider, G., Mielke, A.:
%The validity of modulation equations for extended systems with cubic
%  nonlinearities.
%{Proc. Roy. Soc. Edinburgh Sect. A} {\bf 122}, 85-91 (1992).


\bibitem[LG19]{LG19}
Li, Y., Gao, Y.:
The method of lower and upper solutions for the cantilever beam equations with fully nonlinear terms.
J. Inequal. Appl. {\bf 136} 1-16 (2019). 



%\bibitem[MN13]{MN13}
%Masmoudi, N., Nakanishi, K.:
%Multifrequency NLS scaling for a model equation of gravity-capillary
%  waves.
%{Commun. Pure Appl. Math.} {\bf 66}(8), 1202-1240 (2013).


%\bibitem[P11]{P11}
%Pelinovsky, D. E.: Localization in periodic potentials. From Schr\"odinger operators to the Gross-Pitaevskii equation.
%{London Mathematical Society Lecture
%Note Series} {\bf390}. Cambridge: Cambridge University Press (2011).


%\bibitem[S98a]{Schn98Nodea}
% Schneider, G.: Justification of modulation equations for hyperbolic systems via normal forms.
%{NoDEA Nonlinear Differential Equations Appl.} {\bf 5}, 69-82 (1998). 

%\bibitem[S98b]{S98b}
%Schneider, G.: 
%Approximation of the Korteweg-de Vries equation by the Nonlinear Schr\"odinger equation. 
%{J. Differential Equations}   {\bf 147}, 333-354 (1998). 

\bibitem[S05]{S05}
Schneider, G.:
Justification and failure of the nonlinear Schr\"odinger equation
  in case of non-trivial quadratic resonances.
{J. Differential Equations} {\bf 216}, 354-386 (2005).

%\bibitem[S11a]{Schn11OWbuch}
%Schneider, G.:
%The role of the Nonlinear
%Schr\"odinger equation in
%nonlinear optics. In: 
%Oberwolfach seminars 42:
%Photonic Crystals: Mathematical Analysis and Numerical Approximation
%by D\"orfler, W., Lechleiter, A., Plum, M., Schneider, G. and Wieners, C..
%Birkh\"auser (2011).
%
%\bibitem[S11b]{Schn11} 
%Schneider, G.: Justification of the NLS approximation for the KdV equation using
%the Miura transformation. Advances in Mathematical Physics (2011) 854719.



\bibitem[SSZ15]{SSZ15} 
Schneider, G., Sunny, D.A., Zimmermann, D.:
The NLS approximation makes wrong predictions  for the 
water wave problem  in case of small surface tension and spatially periodic boundary conditions.
{J. Dynam. Differential Equations} {\bf 27}(3), 1077-1099 (2015). 



\bibitem[SU17]{SU17}
Schneider, G., Uecker, H.:
{Nonlinear PDEs.
A dynamical systems approach.}
Graduate Studies in Mathematics {\bf 182}, {American Mathematical Society, Providence, RI},
ISBN: 978-1-4704-3613-1  (2017)


%\bibitem[SW00]{SW00}
%Schneider, G., Wayne, C.E.:
%{ The long wave limit for the water wave problem.
%I. the case of zero surface tension.}
%Comm. Pure Appl. Math. {\bf 53}(12), 1475-1535 (2000). 
%
%\bibitem[SW02]{SW02}
%Schneider, G., Wayne, C.E.:
%{The rigorous approximation of long-wavelength capillary-gravity
%waves.}
%Arch. Rat. Mech. Anal. {\bf 162}, 247-285 (2002). 

\bibitem[SW11]{SW11}
Schneider, G., Wayne, C.E.:
Justification of the NLS approximation for a quasilinear water
  wave model.
 J. Differential Equations {\bf 251}, 238-269 (2011).  

%\bibitem[Sh85]{Sh85}
%Shatah, J.:
%Normal forms and quadratic nonlinear Klein-Gordon equations.
%{Comm. Pure Appl. Math.} {\bf 38}, 685-696 (1985)

%
%\bibitem[T15]{T15}
%Totz, N.:
%{A justification of the modulation approximation to the 3d full water wave problem.} 
%{Comm. Math. Phys.} {\bf 335}(1), 369-443 (2015).

\bibitem[TW12]{TW12} Totz, N., Wu, S.: A rigorous justification of the 
modulation approximation to the 2D full water wave problem. 
{Comm. Math. Phys.} {\bf 310}(3), 817-883 (2012).

\bibitem[WC06]{WC06} Wang, S., Chen, G.: 
Cauchy problem of the generalized double dispersion equation. 
{Nonlinear Anal.} {\bf 64}(1), 159-173 (2006).



\bibitem[Z68]{Z68}
Zakharov, V. E.: Stability of periodic waves of finite amplitude on the surface
of a deep fluid. Journal of Applied Mechanics and Technical Physics, {\bf 9}(2), 190-194 (1968).


\end{thebibliography}

\end{document}